\documentclass{gtpart}     
\title{A Discrete Proof of the General Jordan-Schoenflies Theorem}

\author{Li Chen}
\givenname{Li}
\surname{Chen}
\address{Department of Computer Science and IT, University of the District of Columbia, Washington, DC 20008, USA}
\email{lchen@udc.edu}
\urladdr{www.udc.edu/prof/chen}

\author{Steven G. Krantz}
\givenname{Steven}
\surname{Krantz}
\address{Department of Mathematics, Washington University in St. Louis, St. Louis, MI 63130, USA,}
\email{sk@math.wustl.edu}
\urladdr{http://www.math.wustl.edu/~sk/}

\arxivreference{}
\arxivpassword{}

\volumenumber{}
\issuenumber{}
\publicationyear{}
\papernumber{}
\startpage{}
\endpage{}
\doi{}
\MR{}
\Zbl{}
\received{}
\revised{}
\accepted{}
\published{}
\publishedonline{}
\proposed{}
\seconded{}
\corresponding{}
\editor{}
\version{}

\usepackage{graphics} 
\usepackage{epsfig} 
\usepackage{caption}
\usepackage{float}

\newfam\msbfam
\font\tenmsb=msbm10  scaled \magstep1 \textfont\msbfam=\tenmsb
\font\sevenmsb=msbm7 scaled \magstep1 \scriptfont\msbfam=\sevenmsb
\font\fivemsb=msbm5  scaled \magstep1 \scriptscriptfont\msbfam=\fivemsb
\def\Bbb{\fam\msbfam \tenmsb}

\def\RR{{\Bbb R}}

\def\11{{\Bbb 1}}
\def\KK{{\Bbb K}}

\def\MM{{\Bbb M}}
\def\SS{{\Bbb S}}

\newtheorem{theorem}{Theorem}[section]
\newtheorem{corollary}[theorem]{Corollary}
\newtheorem{proposition}[theorem]{Proposition}
\newtheorem{lemma}[theorem]{Lemma}

\newtheorem{definition}[theorem]{Definition}

 \def\HollowBox #1#2{{\dimen0=#1 \advance\dimen0 by -#2
       \dimen1=#1 \advance\dimen1 by #2
        \vrule height #1 depth #2 width #2
        \vrule height 0pt depth #2 width #1
        \llap{\vrule height #1 depth -\dimen0 width \dimen1}%
       \hskip -#2
       \vrule height #1 depth #2 width #2}}
 \def\BoxOpTwo{\mathord{\HollowBox{6pt}{.4pt}}\;}

\def\endpf{\hfill $\BoxOpTwo$}

\begin{document}

\begin{abstract}
In the early 1960s, Brown and Mazur proved the general Jordan-Schoenflies theorem. This fundamental theorem states: If we embed an $(n-1)$ sphere $S^{(n-1)}$ locally flatly in an $n$ sphere $S^{n}$, then it decomposes $S^{n}$ into two components. In addition, the embedded $S^{(n-1)}$ is the common boundary of the two components and each component is homeomorphic to the $n$-ball.\newline
This paper gives a constructive proof of the theorem using the discrete method.  More specifically, we prove the equivalent statements: Let $M$ be an  $n$-manifold, which is homeomorphic to $S^{n}$. Then, every $(n-1)$-manifold $S$, a submanifold with local flatness in $M$, decomposes the space $M$ into two components where each component is homeomorphic to an $n$-ball. The method was chosen in order to evaluate the computability and computational costs among operations between cells regarding homeomorphism. In addition, methods within the proof can be extended to applications in design algorithms under the assumption that homeomorphic mappings are constructible and computable. In this new revision, We add some new detailed discussions.


\end{abstract}


\maketitle
\newpage
\tableofcontents


\section{Introduction}

The classical Jordan curve theorem was believed to have been first proven by
Veblen in 1905. The Jordan curve theorem states that a simple,
closed curve $C$ separates the plane into two
components. After its first publication, there were many other proofs that followed, including Tutte's proof based on planar graphs in 1978
~\cite{Tut}.  Other discrete proofs have also been explored. In \cite{Che13}, Chen gave a discrete proof of the
classical Jordan curve theorem based on discrete
manifolds that are cell-complexes.

The Jordan-Schoenflies theorem is a further development of
the Jordan curve theorem. The general Jordan-Schoenflies theorem is a fundamental theorem
in geometric topology~\cite{Moi,Bro,Maz}, which states that
embedding an $(n-1)$-sphere $S^{(n-1)}$ locally flatly into an
$n$-sphere $S^{n}$ decomposes the space into two components.
In addition, the embedded $S^{(n-1)}$ is the common boundary
of the two components, where each component is homeomorphic to the
$n$-ball.

To better understand this theorem, we can embed a
$1$-sphere (a circle) into a $2$-sphere (a globe). This could result in two ``bowls'' where the circle is their common boundary. If the circle is a simple closed curve, then it would separate the two components (given that the curve is
locally flat, a situation we later explain).

The Jordan-Schoenflies theorem
confirmed that one component is homeomorphic to an open disk.
However, this theorem is only valid in two dimensional space.
In three-dimensional space, there is a counterexample given by
Alexander's horned sphere: It separates space into two
regions, but the regions are so twisted that they are not
homeomorphic to a normal 3-disk ~\cite{Hat,KS}.

Since the Alexander horned sphere embedding cannot be made
differentiable nor polyhedral, Mazur \cite{Maz,Maz61} used the
concept of ``nice embedding'' to obtain the general
Jordan-Schoenflies theorem. Brown \cite{Bro} simplified
Mazur's concept to local flatness that can prevent the
infinite twists of the Alexander horned sphere.

Compared to the original Jordan curve theorem, the classic
description of the general Jordan-Schoenflies theorem only considers
spheres as the ambient space. An equivalent, more inclusive statement to the theorem is as follows: For an
$n$-manifold $\MM$ that is homeomorphic to an $n$-sphere and
an $(n-1)$-submanifold $\SS\subset \MM$ that is homeomorphic to
an $(n-1)$-sphere, $\SS$ decomposes $\MM$ into two components
where their common boundary is $\SS$ if $\SS$ is locally flat
in $\MM$.

We give a discrete proof of the above description for the general
Jordan-Schoenflies theorem in higher dimensions. We assume
that $M$ is a triangulation or a polygonal decomposition of
$\MM$ and $S$ is an $(n-1)$ dimensional discrete submanifold of
$M$, where $S$ is locally flat, closed, and orientable.

Our proof
will reduce
to
proving the two following theorems: {\bf (1)} (the general
Jordan theorem) Every $(n-1)$-submanifold $S$ that is
homeomorphic to a sphere and is a submanifold with local
flatness in an $n$-manifold $M$, which is homeomorphic to an
$n$-sphere, decomposes the space $M$ into two components.
In other
words, embedding an $(n-1)$-sphere $S^{(n-1)}$ in an
$n$-sphere $S^{n}$ decomposes the space into two components,
and the embedded $S^{(n-1)}$ is their common boundary. {\bf (2)}
Each of the two components is homeomorphic to the $n$-ball.

The advantage of our proof is that we use a completely
constructive method that can also be used to design algorithms for
applications. For instance, we can use this method to actually
deform a separated component into an $n$-ball. This kind of
procedure has potential applications in the massive data
processing of topological structures in persistent analysis.




\section{Concepts Review and New Concepts}

In this section, we review and clarify existing concepts and introduce some new concepts used in the proof of the discrete form of the Jordan-Schoenflies theorem in Section 5.
To begin with definitions, a $k$-cell is basically a $k$-dimensional open manifold
that is homeomorphic to a (an open) $k$-ball. For the purposes of this paper, we add certain constructive properties to the $k$-cell: 1) A $k$-cell and its boundary
must be finite-time constructible, and 2) The $k$-cell and its boundary  are finite-time computable (decidable) in computing science if $k$ is a
fixed number. In other words, even though we do not require
that a $k$-manifold be constructible or computable, we at least require that a $k$-cell and its boundary be constructible and computable.

Constructible, decidable, and polynomial time computable functions are three main categories regarding computability in mathematics~\cite{BB,Sip}.
The word ``computable'' means able to be calculated with a computer program, and ``polynomial time computable'' means able to be calculated with a computer program in polynomial time. Constructible functions may not necessarily be computable.

A cell should be simple enough that it can be constructed or determined effectively.  For example, a simplex (triangle) is a cell that can be
determined easily; it is linear-time (a special case of polynomial-time) computable if its vertices are not located at non-rational number coordinates.
~\footnote{Note that not every real number is computable. For readers who only want to deal with the concept of constructability,
 treat the concept of polynomial time decidable as having a sample equal to "easily constructible."}  A cubic cell with integer coordinates is also easy to determine and is polynomial time decidable.

Therefore, from a computing standpoint, it is not rational to make a cell as complicated as a manifold. If we allow the cell to be very complex, then we may encounter
a difficult argument loop, which is why previous researchers primarily used simplices. However, in this paper, we
use general cells since cubic cells are easy to deal with intuitively. Our definition already includes
the simplex as a special case.

We want to add the polynomial-time computable property of cells to study the structure of manifolds. Therefore, we require a third property of a $k$-cell: 3)
The $k$-cell and its boundary are polynomial-time computable.

This property allows us to use a computing algorithm to determine if a
set is a $k$-cell. We can also find the boundary and other information regarding the $k$-cells in polynomial-time.
In addition, we only consider a finite number of cells.

We can assume that each $k$-cell only contains a finite or constant number of $i$-cells with respect to $k$, $0\le i < k$.
This means that if $k$ is a fixed number, then the number of sub-cells of the $k$-cell will be bounded by a function with respect to
$k$.

\subsection{Triangulation and Simplicial Complexes}

Triangulation is a type of decomposition of a continuous space.
In other words, a 2D or 3D space is partitioned into triangles or tetrahedra, which are both examples of simplices.

Mathematically, an $n$-simplex $\Delta$ is a convex hull of $n + 1$ vertices $v_0,\dots, v_n \in \RR^n$ such that the
vectors $\langle v_1 - v_0 \rangle,\dots, \langle v_n-v_0 \rangle$ are linearly independent (called an affine space).
This simplex is the set determined by

$$
\Delta = \biggl \{\alpha_0 v_0 + \dots+\alpha_n v_n \biggm |  \sum_{i=0}^{n} \alpha_i=1 \mbox{ and } \alpha_i\ge 0 \biggr \}  \mbox{.}
$$

An $(n-1)$-face of $\Delta$ is a subset of $\Delta$ where  $\alpha_k \equiv 0$ for a fixed $k$. The intuitive meaning
is that if $\Delta$ is a 3D  tetrahedron, then a $2$-face of $\Delta$ is a boundary triangle of the tetrahedron.
We can define an  $i$-face of $\Delta$ as the $(n-i)$ numbers $\alpha_j \equiv 0$ in the above equation.

We define (in dimension theory) $\emptyset$ to be a $(-1)$-face of $\Delta$. A simplicial complex is defined as a set of simplices such that:
{\bf (1)}  Any $i$-face of a simplex from $\KK$ is also in $\KK$, and
{\bf (2)} the intersection of any two simplices $\Delta_1, \Delta_2 \in \KK$ is a face of these two simplices in $\KK$.

The dimension of a simplicial complex is defined as the highest dimension of any simplex in $\KK$.  We can see
that a triangulation of a plane forms a 2-simplicial complex. \cite{Lef}


In this subsection, we show that the discrete space used
is a special case of the CW complex when
embedding the discrete space into a Hausdorff space.
In computer graphics or discrete geometry, a 2D complex is represented by three sets: (1) A set of
0-cells, (2) A set of 1-cells, and (3) A set of 2-cells. Each 1-cell is represented by two end points (0-cells),
and each 2-cell is represented by the boundary polygon, which is a closed path formed by several 1-cells.
Therefore, a 2-cell is completely determined by its boundary, which
is a closed 1-cycle.
Intuitively, when we
think about a simple 1-cycle representing a 2-cell,
we can fill substance into the cycle (such as a deformed 2-ball filling the inner part of the cycle). A famous example would be the minimal surface when given a boundary cycle. This process of obtaining this surface need to be computable here.

The following is the formal definition of discrete space:

Let us
define a partial graph $P(S)$, $S\subset V$ to be a subgraph
where each edge $(a,b)$ of $G$ is in $P(S)$ if $a,b\in S$. The
concept of a minimal cycle $C$ is a cycle that does not
contain any proper subset that is also a cycle. Strictly speaking,
we mean that the partial graph of any proper subset of
vertices in the cycle $C$, with respect to the original graph
$G=(V,E)$, does not contain any cycles.

A discrete space is a graph $G$ that has an associated structure. We always assume that $G$ is finite, meaning that $G$ contains only a
finite number of vertices.  Specifically, ${\cal C}_2$ is the set of all minimal cycles representing all possible 2-cells;
$U_2$ is a subset of  ${\cal C}_2$.   Inductively, ${\cal C}_3$  is the set of all minimal 2-cycles made by $U_2$.   $U_3$ is a subset of  ${\cal C}_3$.
Therefore, $\langle G,U_2,U_3,\cdots,U_k \rangle$ is a discrete space, and we can see that a simplicial complex is a discrete space.
For computational purposes, we require
that each element in $U_i$ can be embedded into a Hausdorff space or Euclidean space using a polynomial time algorithm or an efficient constructive method.
Such a mapping will
be a homeomorphism to an $i$-disk with the internal area of the cell corresponding to an $i$-ball that can also be determined in polynomial time.
Another thing we need to point out
is that $u\cap v$ in   $\langle G,U_2,U_3,\cdots,U_k \rangle$ must be connected. In most cases,  $u\cap v$ is either a single $i$-cell in $U_i$ or is empty.

In general,  $u\cap v$ is homeomorphic to an $i$-cell or empty.  In \cite{Che04,Che14}, we use connected and regular points
to define this idea for algorithmic purposes because
homeomorphism is difficult to compute. Now, we require  that $u\cap v$ be
homeomorphic to an  $i$-cell in polynomial computable time.  We also require that deciding whether
an $i$-cycle is a minimal cycle or an $(i+1)$-cell be polynomial time computable as well.
For example, a polyhedral partition can usually be completed in polynomial time in computational geometry.

The CW complex is a special type of cell-complex. Its definition was first introduced by Whitehead in \cite{Whit}.
A more abstract definition of CW complexes
can be found in \cite{Hat}. However, for simplicity, we use Whitehead's original definition.

A cell complex $\KK$ is: {\bf (1)} A Hausdorff space, and {\bf (2)} The union of disjoint (open) cells denoted by $e^{(0)},\cdots, e^{(n)}$
with the following characteristic mapping properties:
Let $e^{(n)}$ be an $n$-cell, meaning that it is homeomorphic to an open $n$-ball $B^n$. Let $D^n$ be the (closed) $n$-disk. We know
that $S^{n-1}=D^n \setminus B^n$ is an $(n-1)$-sphere.
The closure of $e^{(n)}$,  ${\bar e}^{(n)}$, is the image of a mapping $f$ from the $n$-disk $D^n$ to ${\bar e}^{(n)}$ ($f:D^n \to {\bar e}^{(n)}$ ) such that:
{\bf (1)} $f$ is a homeomorphism onto $e^{(n)}$ with restriction to $D^n \setminus S^{n-1}$, {\bf (2)}  ${\bar e}^{(n)} \setminus e^{(n)}$
(denoted by $\partial e^{(n)}$) is a subset of the $(n-1)$-skeleton (or section) of $\KK$.

The $(n-1)$-skeleton (or section) of $\KK$ is usually denoted by $K^{n-1}$, meaning that all cells whose dimension does not exceed $(n-1)$ are in $\KK$.
A CW complex is a special cell complex with
properties called closure-finite and weak topology:  {\bf (1)} For any $e\in \KK$, ${\bar e}$ only intersects a finite number of cells in
$\KK$. (This means that the boundary of $e$ only contains a
finite number of cells in $\KK$). This is called {\it closure-finite}.)  {\bf (2)} A subset $X$ of $\KK$ is closed if and only if
$X\cap {\bar e}$ is closed in $X$ for each cell $e$ in $\KK$. (This is called {\it weak topology}.)

We know that if $\KK$ is finite, then this cell complex is a
CW complex~\cite{Whit}. We would only need to show that a
discrete space, which is always finite as defined above, is a
cell complex. It is easy to show that, in $G=(V,E)$, $V$
contains all 0-cells, and $E$ is the 1-cell set. A simple
1-cycle, which is finite, can be embedded into
Euclidean space (or a Hausdorff space) as the boundary of a
1-ball. If this 1-cycle is in $U_2$ (which must be a minimum
cycle), then it represents a 2-cell with boundaries. The inner
part of this 2-cell is an abstract entity of the cell. It is
represented as an element existing in $U_2$, and it becomes
concrete (or real) when it is embedded into an actual space such as Euclidean
space. In addition, the boundary is made up of 1-cells and
0-cells (in $K^1$, 1-Skeleton).
Such an embedded mapping is a
characteristic map required for cell complexes.

Inductively, for any $e\in U_n$, we know that $e$ is a minimum $(n-1)$-cycle and can be algorithmically
embedded into a  Hausdorff space (as we assume). The inner part of $e$
is homeomorphic to $B^n$ (again we assume this to be algorithmically doable in the construction of $\langle G,U_2,U_3,\cdots,U_k \rangle$).
This $(n-1)$-cycle is a subset of $K^{n-1}=U_0\cup \cdots \cup U_{n-1}$, where  $U_0=V$ and  $U_1=E$.

Triangulations and the piecewise linear decomposition of a space  in Euclidean space are two examples.

The only restriction of the discrete space is that we require that the intersection of the closures of two cells must be homeomorphic to an $i$-cell. For a triangulation, this is true. For a
piecewise linear decomposition, we can usually use an algorithm to refine the original decomposition to satisfy such a property.
This property is somewhat similar to the closure finite
property. The reason we want finiteness in a CW complex is that it is not possible to constructively determine whether a boundary has an infinite number of cells.

The original meaning of an $n$-sphere in the Jordan-Schoenflies theorem is slightly different from the algorithmic homeomorphic
mapping from a simple discrete $n$-cycle. This is because
perfect discrete spheres are hard to describe.
Some similar ideas and historical reviews related to the proofs of this theorem can be found in~\cite{Sie}.



\subsection{Discrete Manifolds}

The concept of discrete manifolds was created for computational purposes. In other words, we re-define a simplicial complex (or a cell complex) as a discrete manifold;
We can view a discrete manifold as a simplicial complex that is made from the decomposition
of a manifold. This simplicial complex only contains a finite number of simplices.

Our definition of discrete space is a special case of piecewise linear space, and the definition of discrete manifolds is more strict. Here, we present a brief description (the formal definition
of discrete space and discrete manifolds can be found in the Appendix in this paper). We can see that a triangle is
 determined by three edges that form a closed cycle. We can say that a 2-simplex is formed by a closed 1-cell path (cycle). This
cycle does not contain any other cycles and is called a minimal cycle. Intuitively, we can fill some materials inside the cycle to make a solid triangle.
However, computationally, there is no need to do this (filling) since we are not going to split a triangle by adding a point or
doing any other surgeries
on a triangle in this paper. In any case, we can say that a 2-cell is determined by a minimal (closed path) cycle of 1-cells.

In general, the boundary of a $k$-cell is a minimal closed $(k-1)$-cycle,  and the $k$-cell is fully determined by its boundary, a minimal closed $(k-1)$-cycle.
However, a minimal closed $(k-1)$-cycle might not be the boundary of a $k$-cell in general discrete space since it is
dependent on whether the inner part of the cell is defined (included) in the complex.  For instance, the inner cycle of a torus is a minimal cycle,
but it is not the boundary of a 2-cell.
(Note that we sometimes simplify by saying that a $k$-cell is a minimal $(k-1)$-cycle because their computational representations are the same.)

To repeat the above idea, we assume that a $1$-cycle is a closed simple path
that is homeomorphic to a $1$-sphere. Further, a $(k-1)$-cycle
is  homeomorphic to a $(k-1)$-sphere. The boundary of a $k$-cell is a $(k-1)$-cycle.
Note that the boundary of a $k$-cell requires homeomorphism to be decided in polynomial time.


We also have further requirements on {\it regular manifolds}. A regular $k$-manifold $M$ must have the following properties:
{\bf (1)} Any two $k$-cells must be $(k-1)$-connected, {\bf (2)} Any $(k-1)$-cell must be contained in one or two
$k$-cells, {\bf (3)} $M$ does not contain any $(k+1)$-cells, and {\bf (4)} For any point $p$ in $M$, the neighborhood of $p$ in $M$,
denoted by $S(p)$, must be $(k-1)$-connected in $S(M)$.

In the theory of intersection homology or piecewise linear topology~ \cite{GM} (or as proven in \cite{Che13}), the neighborhood $S(x)$ of $x$ (which consists of
all cells that contain $x$) is called the {\it star} of $x$. Note that $S(x) \setminus \{x\}$ is called the {\it link}.
Now we have: If $K$ is a piecewise linear (PL) $k$-manifold, then the link $S(x) \setminus \{x\}$ is a piecewise
linear $(k-1)$-sphere.
We can write ${\rm Star}(x)$ as $S(x)$ and ${\rm Link}(x)={\rm Star}(x) \setminus \{x\}$.  In general, we can
define ${\rm Star}({\rm arc}) = \cup_{x\in {\rm arc}} {{\rm Star}(x)}$. Therefore, ${\rm Link}(\rm arc) = {\rm Star}({\rm arc}) \setminus \{arc\}$.
${\rm Star}({\rm arc})$ is the envelope (or a type of closure) of ${\rm arc}$.

We also know that if any $(k-1)$-cell is contained by two $k$-cells in a $k$-manifold $M$, then $M$ is closed.

\subsection{Contraction and Simply Connected Spaces}

A simple path (in a graph) is called a pseudo-curve, or semi-curve in \cite{Che14}. If a pseudo-curve does not
contain all points of a 2-cell, or $k$-cell if $k\ge 2$), then this pseudo-curve is
called a discrete curve. We can see that a discrete curve is similar to the locally flat curve in continuous cases.
Detailed definitions can be found in the Appendix. \footnote{A pseudo-curve is equivalent to a PL curve in
geometry or topology. It is called a pseudo-curve in digital geometry since we can usually collect or sample
discrete points in real world applications. Therefore, for applications, a $k$-manifold is represented by 0-cells. If we want
a unique interpretation of these 0-cells, then we must have some restrictions. In such a case, for example,
we want to eliminate all instances of 2-cells in the
representation of the curve, so we use the concept of the pseudo-curve as the arbitrary curve. We
refer to {\it the discrete curve} as a set of 0-cells not containing all 0-cells of a 2-cell. Otherwise, it is hard to distinguish a curve or a 2-cell computationally. However,
in this paper, without indicating specifically, the pseudo-curve (or discrete curve)
is the PL curve so that there is no confusion when a curve is represented using edges (1-cells) or a $k$-manifold is represented using $k$-cells.
In addition, if all $k$-cells in a complex form the boundary of a $(k+1)$-cell,
then we can view this complex as a closed $k$-manifold if
the $(k+1)$-cell is not included in the complex. We could also view this complex as a $(k+1)$-cell if the $(k+1)$-cell is included in
the complex. We call the former the discrete $k$-manifold and the latter boundary of a $(k+1)$-cell
the discrete pseudo-$k$-manifold.
The reason is that in real world problems, filling a closed $k$-manifold is an abstract matter for computers. This does not affect
the concepts or proofs in this paper.}	

Two (pseudo-)curves $C$ and $C'$ are called {\it gradually varied} if we can deform $C$ to match $C'$ using only
one step. A detailed definition is given in the Appendix.

For the concepts of discrete contraction, we also need the following definition:

\begin{definition} \rm
We say that a collection of simple paths (a pseudo-curve)
is {\it side-gradually varied} if there are no transversal intersections, i.e. no crossovers.
\end{definition}
			
The transversal intersection (see Appendix) is also called  the cross-over in discrete space meaning that
two curves intersect each other, not only touch (each other).


For further calculations, we want to define a special operator ${\rm XORSum}$, which stands for Exclusive-Or-Sum.
${\rm XORSum}$ is $sum(modulo2)$
\ in Newman's book~\cite{New}. $\hbox{\rm XORSum}$ is a computer science term that is relatively easy to understand.
Let $E(C)$ be all of the edges in path $C$.  Then, ${\rm XorSum}(C,C') = (E(C) \setminus E(C'))\cup (E(C') \setminus E(C))$.
 The purpose of this operation is to cut out the shared portion between $C$ and $C'$. The remaining edges will be the collection of cycles when the two
end points of $C$ and $C'$ are the same. If these cycles are boundaries of 2-cells, we can move $C$ to $C'$ in one unit time. In other words,
 $C$ and $C'$ are gradually varied (a discretely continuous move without a jump).

A space is said to be simply connected if any closed  simple path can be deformed to a point on the original curve through a collection of
side-gradually-varied simple paths (pseudo-curves). See the Appendix for details.

(Note that pseudo-curves can be embedded in Euclidean space as simple curves. For our purposes, we use both pseudo-curves and (pure) discrete curves depending on which stage we are at during the process. The difference between them is that pseudo-curves may contain ambiguities when view a curve as a set of points (0-cells). Generally, in discrete space, a pseudo-curve contains all vertices of a 2-cell (or $k$-cell) but not necessarily the 2-cell, but for the pure discrete curve, if it contains all vertices of a 2-cell, it contains everything within these vertices including the
2-cell.

However, when we involve higher dimensional cells in a cell complex, we can use appropriate involvements of 2-cells or 3-cells to limit the number of possible outcomes we encounter for pseudo-curves. Therefore, using pseudo-curves for our purposes would not be a problem since we treat pseudo-curves the same as pure curves, so this does not affect our results or proofs.
However, we will need to add a task to select cells in a cell complex.  Pure curves will be simpler.
To be more specific, when we select a curve in the beginning or produce a final curve, we require that the curves be pure discrete curves. However, during the middle of the process, for example during a contraction, we can use pseudo-curves---this is because we can easily determine which 2-cell to be excluded.

As long as we embed a curve in Euclidean space or use higher dimensional cells in a cell complex, the distinction between pseudo and pure curves are nonexistent. Only in discrete space, we want to choose a pure, locally flat curve (not a pseudo-curve) at the beginning to minimize preprocessing, such as memory use. Making a pseudo-curve locally flat would require even more preprocessing time, so we might as well start with the ideal curve at the beginning. However, during the process, allowing pseudo-curves would simplify the method, especially in contractions.)

\subsection{Graph-Distances and Cell-Distances}

In a graph, we refer to the distance as the length of the shortest path between two vertices.
The concept of {\it graph-distance}, in this paper, is the edge distance, meaning how many edges are needed to get from one vertex to another.
We usually use the length of the
shortest path between two vertices to represent the distance in graphs. In order to distinguish from distance in Euclidean space,
we use graph-distance to represent lengths in graphs for this paper.
							
Therefore, graph-distance is the same as edge-distance or 1-cell-distance, which is the number of 1-cells that are needed to travel from one point (vertex) $x$ to point $y$.  We can generalize
this idea to define 2-cell-distance by counting how many 2-cells are needed to go from $x$ to $y$.  In other words,
2-cell-distance is the length of the shortest path of 2-cells that contains $x$ and $y$. In this path, each adjacent pair of
2-cells shares one 1-cell, which means that the path is 1-connected (meaning that two adjacent elements in the path share a 1-cell).

We define $d^{(k)}(x,y)$, the $k$-cell-distance from $x$ to $y$, as the length of the shortest path (or the
minimum number of $k$-cells in such a sequence),
where each adjacent pair of two $k$-cells shares a $(k-1)$-cell. This path is $(k-1)$-connected; in other words, two adjacent elements in the path share a $(k-1)$-cell.

We can see that $d^{(1)}(x,y)$ is the edge-distance or graph-distance, and we have $d(x,y) = d^{(1)}(x,y)$.

(We can also define $d^{(k)}_i(x,y)$) to be a $k$-cell path that is $i$-connected.
However, we do not need to use such a concept in this paper. )

\subsection{Local Flatness}

The concept of local flatness for embedded submanifolds is similar to the smoothness of
manifolds, but it is a stronger definition in some sense. The concept is as follows:

Suppose a $k$-manifold $M_k$ is embedded into an $n$-manifold $M_n$, $k<n$.    The manifold $M_k$ is said to be
{\it locally flat} at $x\in M_k$ if there is a neighborhood  $U_{x} \subset M_n$ of $x$
such that the topological pair $(U_x, U_x\cap M_k)$ is homeomorphic to the pair $( {\mathbb R}^n,{\mathbb R}^k)$,
with a standard inclusion of ${\mathbb R}^k$ as a subspace of ${\mathbb R}^n$.
That is, there exists a homeomorphic mapping $f: U_x \to {\mathbb R}^n$ such that  $f(U_x \cap M_k) = {\mathbb R}^k$.
Here, ${\mathbb R}^n$ is the standard $n$-dimensional real vector space.
				
Let us assume that $M_n$ is closed (since $x$ should not be selected as a boundary point of $M_n$). If every point  $x\in M_k$ is
locally flat in $M_n$, then $M_k$ is called locally flat in $M_n$.

\begin{proposition}\label{PropBrown}[Brown (1962)]  \sl
 If  $k = n-1$, then a locally flat $M_{k}$ is {\it collared}, meaning that it has a neighborhood, which is homeomorphic
 to $M_k \times [0,1]$ where  $M_k$ is (homeomorphically) corresponding to $M_k \times \left \{ \frac{1}{2} \right \}$.
\end{proposition}

This result by Brown provided an intuitive interpretation of local flatness. The meaning of collared is shown in Fig.\ 1 .
The boundary of a collar cannot intersect itself, a property that will be used in the discrete case.

\begin{figure}[h]
	\begin{center}

   \epsfxsize=3in
   \epsfbox{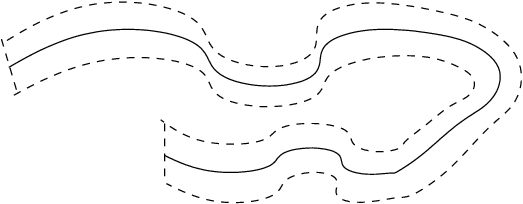}   

\caption{An example of a collar of a curve in continuous space. The boundary of the collar will not intersect itself.}

\end{center}
\end{figure}

In the next section, we discuss the concepts of local flatness and collars for the discrete case.


\section{Local Flatness of Manifolds in Discrete Space}

We have explained the intuitive meaning of a collar of a continuous curve: there
is a neighborhood (of the curve) that does not
intersect itself. The boundary of the collar has a unique distance to
the curve, called graph distance or cell distance.
In this section, we define the local flatness of a
curve or a manifold in the discrete case.

A discrete curve or a curve in discrete space is usually represented as a path of vertices where
each adjacent pair forms an edge. See Section 2. The collar of a discrete curve is a 2-dimensional manifold.
In Fig.\ 2 (a), we present a collar of a discrete curve, but in Fig.\ 2 (b), the boundary of the collar
intersects at a vertex, which means that the latter is not a valid collar.

\begin{figure}[h]
        \begin{center}

   \epsfxsize=5.5in
   \epsfbox{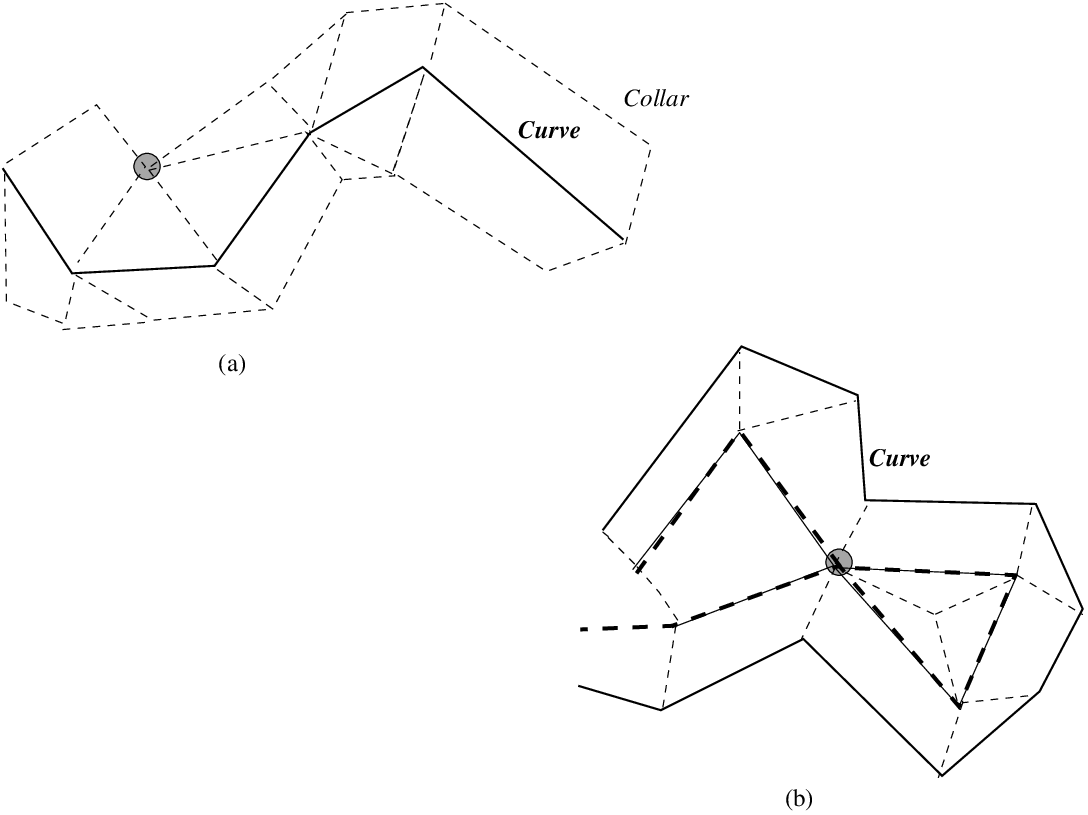}   

\caption{(a) An example of a collar of a curve in discrete space. (b)
The boundary intersects itself, which makes it not locally flat and not a collar.}
\end{center}
\end{figure}

\subsection{Some Observations on Local Flatness in Discrete Space}

Before we give the formal definition of local flatness in discrete
space, we first give an example to illustrate
properties that local flatness should inherently hold.
In Fig.\ 3, points $A$ and $B$ are linked by an edge so that any neighborhood of
$A$ will contain $B$. That is to say that curve $\cdots ACB\cdots$ is not
a locally flat curve if we do not allow a vertex to be added on the
edge $(A,B)$.

\begin{figure}[h]
        \begin{center}

   \epsfxsize=5in
   \epsfbox{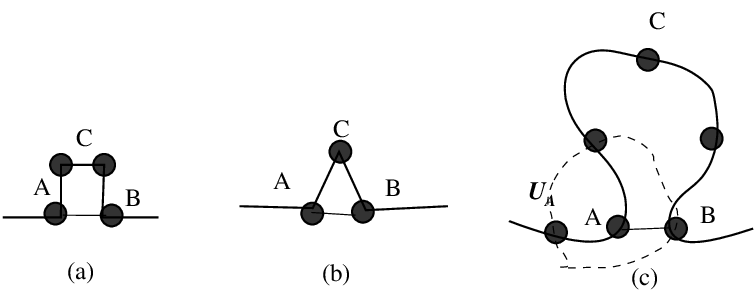}   

\caption{Examples of not locally flat curves: (a) The curve contains points $AC$ and $B$, (b) The curve contains points $ACB$, and (c)
    The curve is of the form $\cdots A\cdots C\cdots B$. }

\end{center}
\end{figure}

Mathematically, this is related to
the definition of local flatness in Section 2.6. We require
$U_x$ to be topologically equivalent (or homeomorphic) to $\RR^n$. We also require

$U_x\cap M_k$ to be topologically equivalent to $\RR^k$.

In the discrete case, $U_A\cap M_k$   is a 2-cell not topologically equivalent to $\RR^1$ in
Fig.\ 3 (a) and (b). (Based on the definition of
discrete manifolds in Section 2, we can further extrapolate that when a set contains all points
of a cell, then this set will also contain the cell. The exception will be in the cases of contraction motions or process. That is also the main difference between a pure discrete manifold and the pseudo discrete manifold. )

On the other hand, in Fig.\ 3 (c),
$U_A\cap M_k$ is two 1-cells plus a single point (0-cell) $B$ that is not
homeomorphic to $R^1$. Here, $M_k$ is a curve.

Therefore, $U_A\cap M_k$ should contain the edge $(A,B)$. That is to say that $U_A\cap M_k$ is a tree where point $A$ connects to three
points, including $B$.

Thus, separating two points on a curve with a vertex not on the curve is one way to maintain the local flatness property.
Topologically, we can see that for some discrete cases, only having one vertex between the two points on the curve is not sufficient. See Fig.\ 4. $U_A\cap M_k$ will intersect $U_A\cap M_k$ at point $C$.

\begin{figure}[h]
        \begin{center}

   \epsfxsize=3in
   \epsfbox{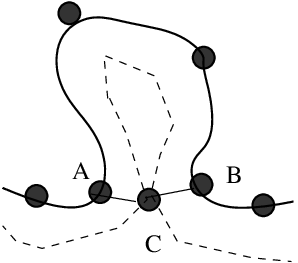}

\caption{A case that satisfies the original definition of local flatness in continuous space but may cause conflicts in discrete space.}

\end{center}
\end{figure}

This is because topologies in discrete space are finite topologies. There is no clear distinction between
the open set and the closed set. When a point (vertex) is contained within two or more collars of a curve,
it may create a more difficult case. This case would not occur in continuous space since
we could use open sets that would not contain such a midpoint $C$ (or a midline). However, we would not be able to store all the points of an open set in finite space (this is not constructive). The open set is a type of imaginary interpretation for continuous space
from the discrete space point of view.


\subsection{Formal Definitions of Local Flatness}

In this subsection, we give the formal definition of local flatness in
discrete space. We use a curve (1D-manifold) as an example for discussion.
Then, we extend the definition to more general cases.

According to the discussion above, an intuitive definition of local flatness is that each pair of
points, if they are not adjacent, in curve $C$ will have a
graph-distance of 3. However, the following example shows that the corner of a curve can be
locally flat. See Fig.\ 5.  Therefore, the challenge is that sometimes one point separation is enough.

\begin{figure}[h]
        \begin{center}

   \epsfxsize=3in
   \epsfbox{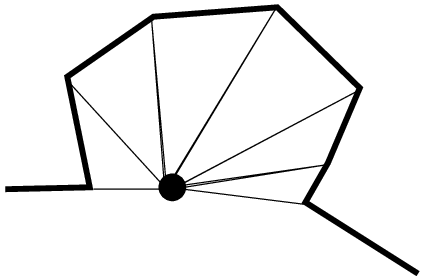} 

\caption{Example of a curve that can be locally flat with one point in between.   }

\end{center}
\end{figure}

This example shows us that
the graph-distance of 3 does not apply to some corner points for local flatness.  In other words,
we could have a collar.

However, for similar cases (see Fig.\ 6 (a)) we would not have a collar since the collar must intersect at a point.
Our intention is to explain why in the following discussion.

\begin{figure}[h]
	\begin{center}

   \epsfxsize=5in
   \epsfbox{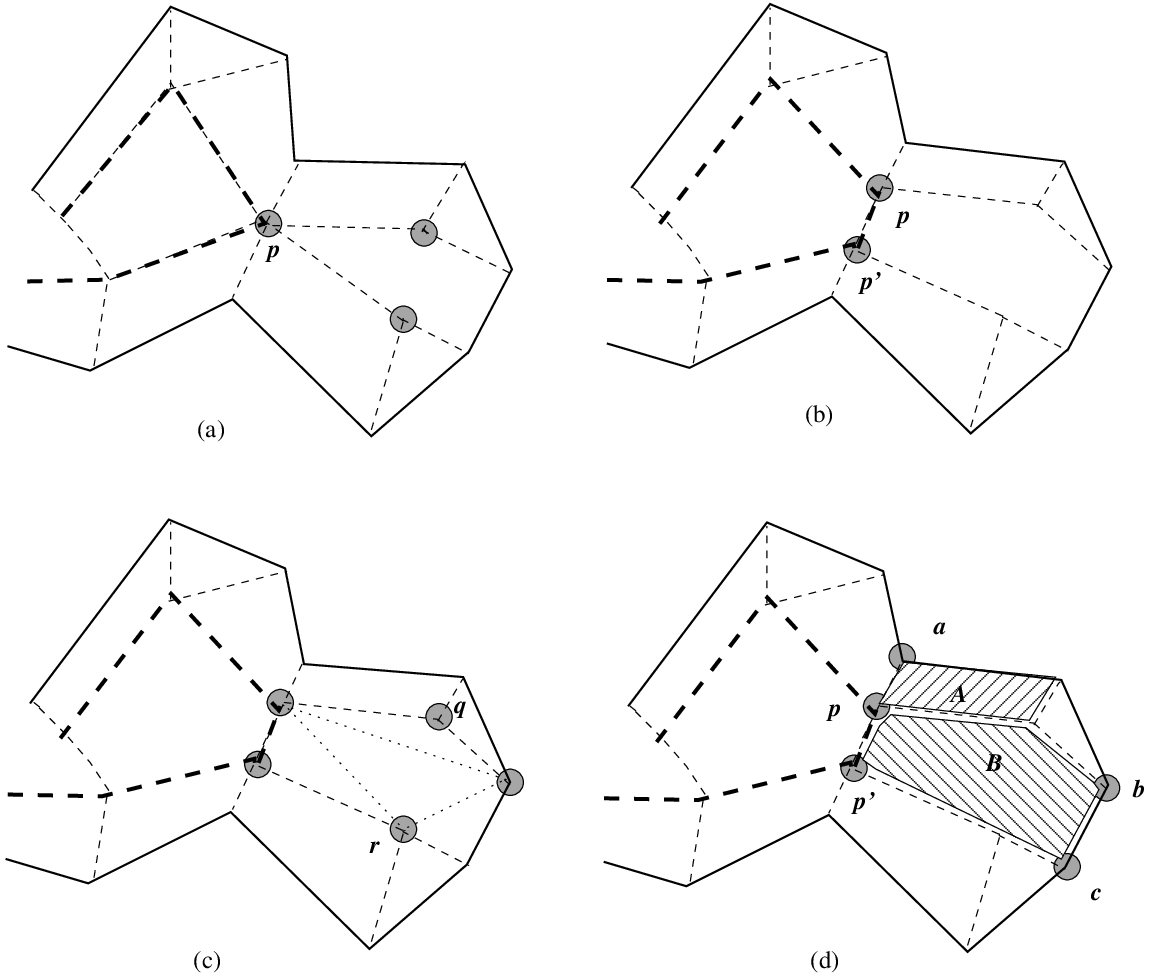}

\caption{Example of a curve that does not have a collar. (a) The collar intersects at a point $p$.  (b)
A point is inserted, but it still does not have a collar.  (c) Triangulation that might not work as well
 unless there is a link between two certain points $q$ and $r$. (d) The reason behind is that the 2-cell-distance between these two points is still 2.}

\end{center}
\end{figure}

We know that Fig.\ 6 (a) does not have a collar since the collar intersects at a point $p$.
When we insert a point $p'$ in Fig.\ 6 (b), we have met the previous
observation of having a graph-distance of 3 for the edges. However, we still could not get a collar.
In Fig.\ 6 (c), we make a triangulation that might not work as well
unless there is a link between two certain points $q$ and $r$.

Fig.\ 6(d) shows that the 2-cell distance is still 2 in this figure for
two points $a$ and $c$ that are not adjacent on the curve. This means that, even if the edge-distance (or graph-distance)
is 3, the two points can still be linked by two 2-cells.

Now, we can conclude that the shape of the 2-cells (in "holding" space) are not the problem. The key is that
a distance of 3 is required for both 1-cell distance (graph-distance) and 2-cell distance.

It is interesting to note that in a triangulated decomposition of a plane, this problem would not
exist. We prove later that in triangulated manifolds, the problem related to 2-cell
distance disappears.

To summarize mathematically: (1) For each pair of two points $p$ and $q$ in $C$, if $p$
and $q$ are adjacent in $M$, i.e. $d_{M}(p,q)=1$, then $p$ and $q$ must be adjacent in $C$. (2) If
$d_{M}(p,q)=2$,
then there must be a point $a$ in $M$  such that $\hbox{\rm Link}(a)\cap C$ is an arc
that contains $p$ and $q$. (Point $a$ is called a focal point of $C$, and focal
points are always at the collar boundary.)
(3) If $d_{M}(p,q)=3$, then for triangulated manifolds, there is a
collar. For general shape 2-cells, we require a distance of 3 in 1D and
2D distance.

In general, for the local flatness of a curve $C\in M$, we require that the 1-cell distance is 3, 2-cell distance is 3, and $k$-cell
distance is also 3 for all cases where $M$ is a $k$-manifold.

The following is a formal definition for local flatness. As we know, a
triangulated discrete manifold $M$ is a discrete manifold where
each face is a simplex.

\begin{definition} (Triangulated discrete manifolds)
A curve $C$ is locally flat in a 2D triangulated discrete manifold $M$
(or $C$ can be a $(k-1)$-manifold in a $k$-triangulated manifold $M$.)
if every pair of points $p$ and $q$ on the curve satisfies one
of the following conditions: (1) $d_{M}(p,q)\ge 3$, (2)
$d_{M}(p,q)=2$, if every point $a$ in $M \setminus C$ on any path with length 2 (from $p$ to $q$)
satisfies $\hbox{\rm Link}(a)\cap C$ is an arc
containing $p$ and $q$, or (3) $d_{M}(p,q)= 1$ where $p$ and $q$ are
adjacent points in $C$, i.e.   $d_{M}(p,q)= d_{C}(p,q)=1$.
\end{definition}

In Definition 3.1, condition (2) describes a point $a$ that is close to
a corner of this curve. If $path_{C}(p,q)$ in $C$ denotes the path from $p$ to $q$ and if  $path(p,q)$ surrounds $a$ (meaning that $\hbox{\rm Link}(a)\cap C=path(p,q)$), then
$d_{M}(p,q)= 2$ is allowed. We can call $a$ a {\it near-corner point with respect to $C$} or a ``focal point'' of the corner of the curve. Any point with such a property must also be near the corner of the curve.
In Fig.\ 7 (a), point $a$ satisfies condition (2), but point $b$ does not. Therefore, it is not a locally flat curve.

Intuitively, if a shortest path passes a point $b$ in $M \setminus C$ that is not a {\it near-corner point with respect to $C$}, then this path is referred to as the ``waist'' of the curve. We
demand that any waist have at least a length of 3.

\begin{figure}[h]
\begin{center}

   \epsfxsize=4.in
   \epsfbox{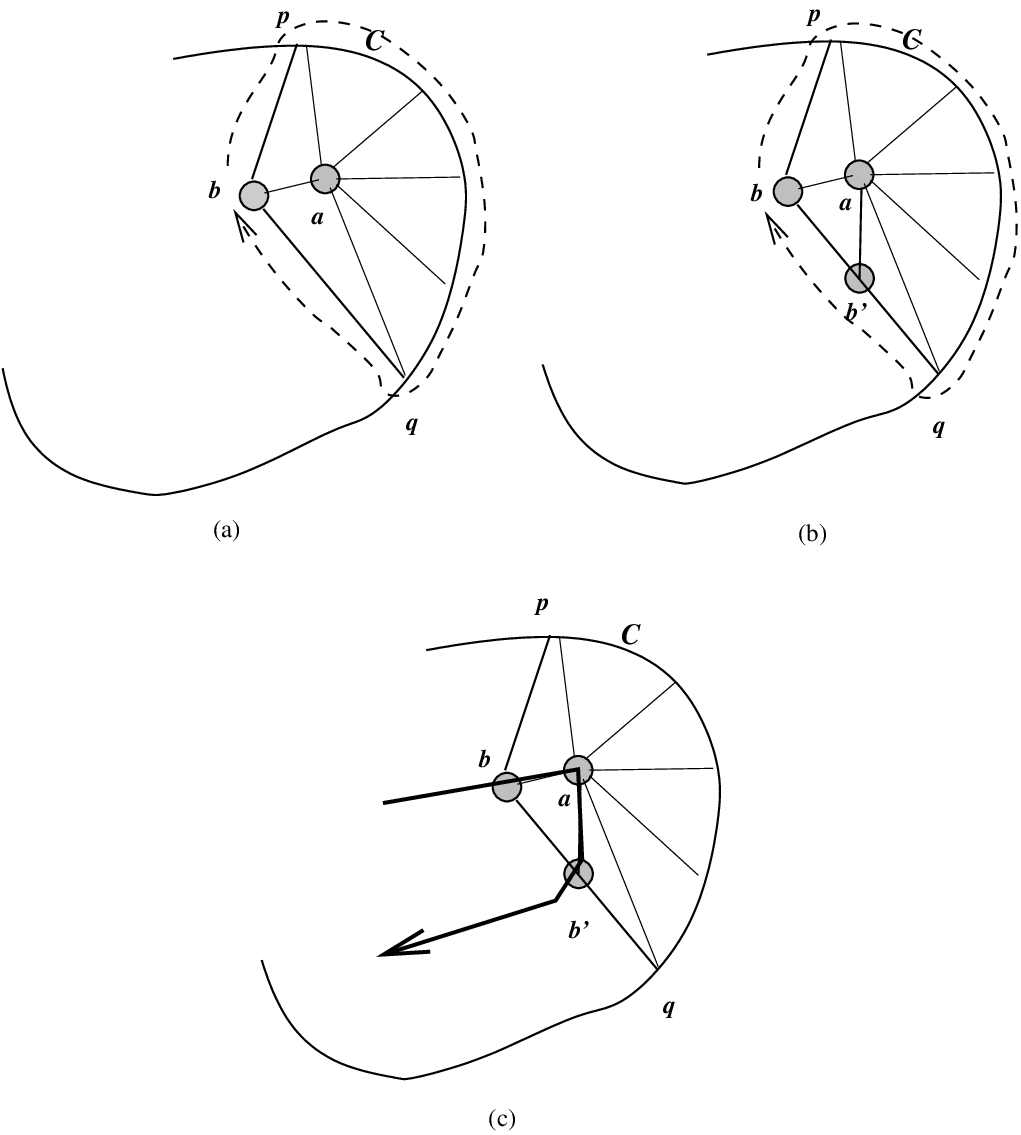}

\caption{Collar and gradually varied properties: (a) Point $p$ does not satisfy Definition 3.1 and $a$ cannot be a collar point, (b) Point $a$ as a collar point, and (c) A collar line based on (b).}

\end{center}
\end{figure}

\begin{definition} (General discrete manifolds)
A curve or submanifold $C$ is locally flat in a discrete $k$-manifold
$M$ if every pair of points $p$ and $q$ on the curve satisfies
one of the following conditions: (1) $d_{M}^{(i)}(p,q)\ge 3$ for all
$i\le k$, (2) $d_{M}^{(i)}(p,q)=2$ for some $i\le k$  if every point $a$ in $M \setminus C$ on such a path satisfies that
$\hbox{\rm Link}(a)\cap C$ is an arc containing $p$ and $q$, or (3)
$d_{M}(p,q)= 1$ where $p$ and $q$ are adjacent points in $C$, i.e.
$d_{M}(p,q)= d_{C}(p,q)=1$.
\end{definition}

Note that the dimension of $C$ in Definition 3.2 above is smaller than $k$.

\begin{lemma}
For a triangulated manifold $M$, graph-distance (1-cell distance) is the same as
$k$-cell distance in $M$, which is a $k$-discrete manifold.
\end{lemma}

\begin{proof}
Proving this is not hard using mathematical induction.
We know that every pair of vertices in a simplex is adjacent.

If $d_{M}(x,y)=d_{M}{(1)}(x,y)=1$, then this means that there is an edge linking $x$ and $y$.
It is obvious that $x$ and $y$ are in a simplex in $M$. Therefore, $d_{M}{(i)}(x,y)=1$ for all $i$,
$1\le i\le k$.

Let us assume that $d_{M}(x,y)=d_{M}{(i)}(x,y)$ is valid for distance $t$. We want to prove that
when $d_{M}(x,y)=t+1$, then $d_{M}{(k)}(x,y)=t+1$.  Since $d_{M}(x,y)=t+1$, there must be a shortest
path $x=x_0,x_1,\cdots, x_t, x_{t+1}=y$ in $M$, so $x_t, x_{t+1}$ must be in a simplex $K$.
Therefore,  $d_{M}{(k)}(x,y)\le t+1$. If $d_{M}{(k)}(x,y)= t$, then $x_{t+1}$ must be contained
in $t$ $k$-simplices that also contain $x(=x_0),x_1,\cdots, x_t$.

According to the assumption, $d_{M}(x,x_t)=d_{M}{(k)}(x,x_t)$  and every pair of points in a simplex
must have an edge linking this pair. The point $x_{t+1}=y$ must be contained in a simplex that also contains three vertices
in the path  $x=x_0,x_1,\cdots, x_t, x_{t+1}=y$. Then, $d_{M}(x,y)=t$ and is not equal to $t+1$.
This gives a contradiction.

It is not difficult to see that $d_{M}{(i-1)}(x,y)\ge d_{M}{(i)}(x,y)$. Therefore, for all $i$, we have
$d_{M}(x,y)=d_{M}{(i)}(x,y)=t+1$.

\end{proof}

Now,  we repeat the formal definition of collar in discrete space:

\begin{definition}
 The meaning of collar in discrete
space is that each point on the collar boundary is 1-cell distance from a point on
the curve $C$. The boundary of the collar does not intersect itself.
\end{definition}

The following lemma is a natural consequence of the above definition.
\begin{lemma}
For an open curve $C$, if it has a collar, the boundary of the collar
must consist of two curves, each of which
is gradually varied to $C$.
However, no point on the collar is
in $C$.
\end{lemma}


In the next subsection, we prove the following theorem:
If $C_{k-1}$ is a discrete local flat $(k-1)$-submanifold in a $k$-discrete manifold $M$, then there is a collar for $C_{k-1}$.


\subsection{Properties of Discrete Local Flatness}

Intuitively, if a discrete curve has a collar, then it is locally flat. However, this observation might not be true
for continuous space.  In this section, we prove that
if a curve is locally flat, then there is a discrete collar.

According to the definition of the discrete collar,
  we consider the boundary of the union of $\hbox{\rm Link}(p)$, $p\in C$. We can see that the boundary of the
collar of $C$ is the subset of $\cup_{p\in C} \hbox{\rm Link}(p)$.

\begin{theorem} (Discrete version of the collar theorem)
       If a curve $C$ is discrete locally flat in a discrete 2-manifold $M$, then there is a collar for $C$.
\end{theorem}

\begin{proof}

We examine two adjacent points $x,x'$ in $C=C(a,b)$, where $a$ and $b$ are two end points.
Here, $C$ is a locally flat arc, meaning that $d_{M}(a,b) \geq 3$ if $a$ and $b$ in $C$ are not adjacent.
The exception is when $a$ and $b$ are in $Star(p)$, $p\in M \setminus C$,
as stated in Definition 3.2.

Let $(x,x')$ in $C$ be a 1-cell. Then, $\hbox{\rm Link}((x,x'))$ is a simple cycle, which is
the union of $\hbox{\rm Link}(x)$ and $\hbox{\rm Link}(x')$.  If $x'$ is not an end point of $C$, then we have another $y\in C$ that is adjacent to $x'$.
(We assume that $C$ is a closed curve, which is closer to the original requirements in Brown's theorem on collars.)

We only want to prove that there are two local curves $B_1$ and $B_2$ that do not intersect $C$ and that the maximum distance (cell-distance) of
each point in $B_1$ or $B_2$ to $C$ is 1. Therefore, $x$ and $y$ have a distance of 2 in $C$ because, if $d_{M}(x,y)=1$, then $C$ would not be locally flat by definition.

There are two cases:
(1) If there is a point $p\in M \setminus C$ such that $\hbox{\rm Link}(p)$ contains $x$ and $y$,
then $\hbox{\rm Link}(p)$ would contain $x'$ and $p$ would be on one side of a ``collar''
($B_1$ or $B_2$) locally. We can always find a new $p$ such that
$\hbox{\rm Link}(p)$ contains $x$,$y$ and $x'$.

This is because if $\hbox{\rm Link}(p)$ does not contain $x'$ and since $\hbox{\rm Link}(p)$ is a simple cycle, then we will have another path from $x$ to $y$ in
$\hbox{\rm Link}(p)$ due to the following: (a) If $(x,y)$ is an edge in $\hbox{\rm Link}(p)$, then $C$ is not locally flat. (b)
This path, denoted by $\rho$, and $x,x',y$ will form a cycle, and now
$\hbox{\rm Link}(x')$ will contain a point in the path.
Using this point as $p$, we can
 use $\hbox{\rm link}(x')$ to select a point not in $C$ to be the new $p$. If we continue doing this, we
can always get a $p$ such that
$\hbox{\rm Link}(p)$ contains $x$,$y$ and $x'$.

This also follows the definition of local flatness that $x$ and $y$ must be connected in $\hbox{\rm Link}(p)$, which means
that $\hbox{\rm Link}(p)$ must also contain $x'$.  Therefore, $x x' y$ is a corner point.

(2) If there is no such point $p$ such that $\hbox{\rm Link}(p)$ contains $x$ and $y$, then
the path from $x$ to $y$ without passing $x'$ will be at least a distance of 3 based
on Definition 3.2. There is a nonempty intersection
of $\hbox{\rm Link}((x,x')) \cap \hbox{\rm Link}((x',y))$, which is a subset of  $\hbox{\rm Link}(x')$. This subset
is a connected curve that does not contain any point in $C$.
This set is the collar on both sides of $x x' y$, $B_1$, and $B_2$.

Now, we use mathematical induction for the rest of proof.
(It is similar to the proof in paper \cite{Che13}. )
If there are $x_1,\dots,x_{k-1}$ satisfying the condition of
having a collar on both sides, then we add an $x_k$ adjacent to $x_{k-1}$.
If there is a point $p$,
$p\in M \setminus C$, such that $x_{k-2}$, $x_{k-1}$, and $x_k$ are neighbors of $p$,
then this $p$ is a corner point. Therefore, if $p$ is a collar boundary point for $x_{k}$,
then the other boundary points are in $\hbox{\rm Link}(x_k)$. In particular, starting
at $p$, all points
with a common 2-cell $x_{k}$ will be in the new collar.

If there is no such $p$, then every cell-distance to points
(other than $x_{k-1}$) will be 3 or larger.
There must be two points on each side of $C$
that are intersection points of $\hbox{\rm Link}(a=x_0,..., x_{k-2})$
Otherwise,
$\hbox{\rm Link}(a=x_0,..., x_{k-2})$  will not reach $\hbox{\rm Link}(x_k)$, and
the distance of $x_k$ to some $x_i$ will be 2.
Therefore, $\hbox{\rm Link}(x_{k-1})\cup \hbox{\rm Link}(x_k)= \hbox{\rm Link}((x_{k-1},x_k))$
(the link of each $n$-cell is a cycle in the cell complex) must have
two points $q_1$ and $q_2$  in $\hbox{\rm Link}(x_{k-1})\cap \hbox{\rm Link}(x_k)$
such that $q_1, y_0,\dots,y_t, q_2$ is an arc (connected) in $\hbox{\rm Link}(x_k)$
but $y_0,...y_t$ are not in $\hbox{\rm Link}(x_{k-1})$.

We have proven that local flatness implies a discrete collar
for a closed curve.

\end{proof}

In general,

\begin{corollary}
Local flatness implies the existence of a discrete collar for each dimension.
\end{corollary}

\begin{proof}
This proof is a continuation of Theorem 3.6. We want to prove the following:
If $C_{k-1}$ is a discrete local flat $(k-1)$-submanifold in $k$-discrete manifold $M$, then there exists a collar for $C_{k-1}$.

Using the same principle as in the above proof of  Theorem 3.6, we can prove this corollary.  Let $C_{k-1}$ be a discrete local flat $(k-1)$-submanifold.
(Precisely, $C_{k-1}$ is assumed to be a discrete pseudo-submanifold meaning that $C_{k-1}$ might contain all 0-cells of a $k$-cell. See the Appendix for details, which will not affect our proof.)
We can assume that  $C_{k-1}$ is closed, i.e. $C_{k-1}$ is a ${(k-1)}$-cycle. We consider two sets: (a) For each
$(k-1)$-cell $Ce\in C_{k-1}$, we have $Link(Ce)$, and (b) For each point $a$ in $M-C_{k-1}$, we consider such an $a$ so that there is a
$k$-cell ${\cal K}\in M$ containing both $a$ and some point $p\in Ce$ where $Ce$ is a ${(k-1)}$-cell in $C_{k-1}$. The intuitive meaning of the set containing all $a$'s is
the neighboring set $C_{k-1}$. In other words, $a$ is a neighbor of the ${(k-1)}$-cycle $C_{k-1}$. For convenience, we use $C$ to represent $C_{k-1}$ in the rest of the proof.

We can prove that all points $a$'s (0-cells), denoted by set $B$ where each element of $B$ has a $k$-cell distance of 1 to $C$ ($C_{k-1}$), form two $(k-1)$-submanifolds.
(i) Assume that $a\in B$. Then, $Link(a)-C$ must contain at least two points. Otherwise, if there is only one element $b\in Link(a)-C$, there must be
a path from one point $p$ in $C$ to another point $q$ in $C$ that passes $b$ on $Link(a)$ . Since $Link(a)$ is a $(k-1)$-cycle, $pbq$ has a cell-distance of 2.
Based on Definition 3.3,  $C$ is not locally flat.  On the other hand, there are at least two points in $Link(a)-C$ that have a cell-distance of 1 to $C$.
Therefore, a part of $B$ is a $k$-cell-connected set on the same side as $a$.
(ii) Since $C$ includes a $(k-1)$-cell containing $p$, this cell is contained in two $k$-cells in $M$, one that also contains $a$ and another that contains a point $e$. Therefore,
$B$ contains a part that is a $k$-cell-connected set on the same side as $e$.
(iii) We want to prove that it is impossible for a point $b\in B$  to be cell-connected to $a$ or $e$.  This is because $Link(x)$ and $x\in B$ must intersect with another
 $Link(y)$ for some $y\in B$. It is obvious that $x$ and $y$ are cell-connected by passing $C$. For instance, we have a $(k-1)$-cell-sequence $Ce_1,\cdots, Ce_t$ fully contained in $C$ where $x$ is in $Link(Ce_1)$
and $y$ is in $Link(Ce_t)$.  $Link(Ce_i) \cup Link(Ce_{i+1})-Link(Ce_i) \cap Link(Ce_{i+1})$ is similar to the connected-sum of two  $(k-1)$-cycles. Therefore, $x$ can reach $y$ by passing the points on $Link(Ce_{j})$,
$j=1,\dots,t$,
without containing any points in $C$. The rest of the details of the proof are similar to the steps in the proof for Theorem 3.6.

\end{proof}

We can easily see that $B_1$ and $B_2$ in the proof of Theorem 3.6 are two gradually varied
curves of $C$. When $C$ is an open curve, $B_1$ and $B_2$ are
connected. (We do not have to use the property of collars to prove
the main theorem of this paper in Section 5. The concept of local flatness is
enough.)

Since it is not difficult to know whether there is a collar (and whether it is locally flat) in (finite) discrete space,
we have the following:

\begin{lemma}
 The necessary and sufficient condition for the existence of
 a discrete collar of a curve $C$ is that it is locally flat.
\end{lemma}

Therefore, the collar condition in the discrete case is the same as local flatness, which is
stronger than the definition of the collar condition in the continuous case.


\section{Properties of the Gradual Variation of Discrete Manifolds}

Deformation is the continuous motion from one object to another.
In mathematics, we can view an object as a manifold.
A special type of deformation, called a contraction, is a cycle on a surface that
shrinks continuously to a point.

Gradual variation is a discrete term for continuity\cite{Che12}, which we define in the Appendix.
Two discrete curves are gradually varied if we can obtain the second curve by moving
each point on the first curve by a distance of at most one unit. In other words,
if $C$ and $C'$ are two discrete curves, we can change $C$ into $C'$ by
moving each point of $C$ by at most a distance of 1.  See more in the Appendix of this paper.

The following is an example to explain this concept.


\begin{figure}[h]
	\begin{center}

   \epsfxsize=4.in
    \epsfbox{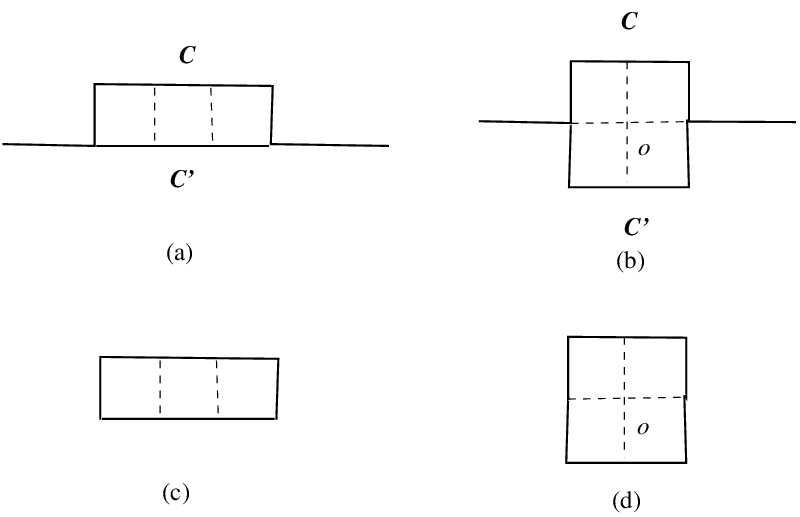} 

\caption{Gradually varied moves of curves: (a) $C$ and $C'$ are gradually varied,
(b) $C$ and $C'$ are not gradually varied, (c)  ${\rm XorSum}(C,C')$ contains all points in three 2-cells, and (d) ${\rm XorSum}(C,C')$
is a cycle that does not contain any 2-cells. }
\end{center}
\end{figure}

The gradual variation between two (discrete) curves can also be described as the modulo-2 sum of two curves (see Appendix).
The modulo-2 sum has another name called the exclusive-OR operation, represented by ${\rm XorSum}(C,C')$.
${\rm XorSum}(C,C')$ contains all edges of $C \setminus C'$ and $C' \setminus C$.
In Fig.\ 8 (c), ${\rm XorSum}(C,C')$ contains all points in three 2-cells.
As we discussed in Section 2 and in \cite{Che04,Che14}, we say that ${\rm XorSum}(C,C')$  also contains these 2-cells.
Fig.\ 8 (d) only contains some edges.  Therefore, two curves are gradually varied if and only if the ${\rm XorSum}$ of them is the union of
several 2-cells (see Appendix).

Discrete deformation can be done by a sequence of gradually varied moves.
In this section, we give some useful properties of the
gradual variation of discrete manifolds.

\subsection{Minimal Gradual Variation Between Curves }

Let us continue the example shown in Fig.\ 8 (a). A new observation can be made in Fig.\ 9.

\begin{figure}[h]
	\begin{center}

   \epsfxsize=4.in
    \epsfbox{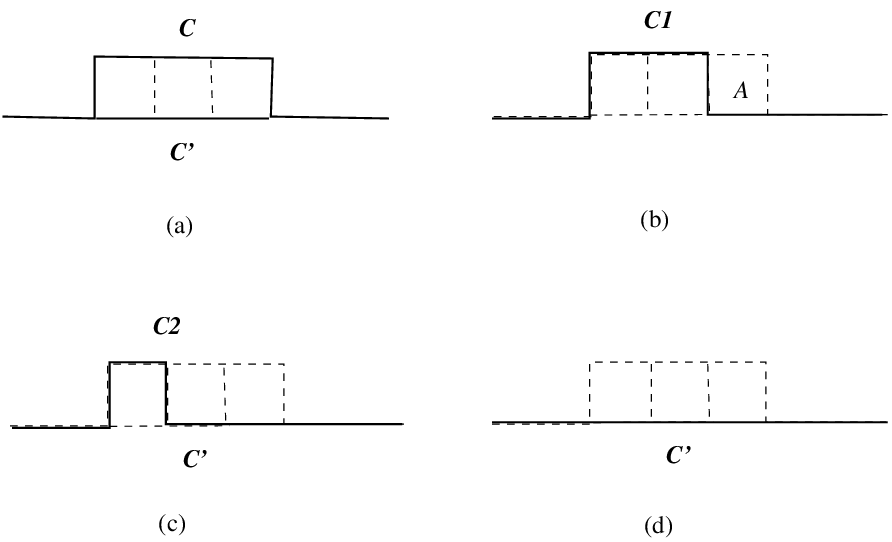} 

\caption{More details of gradually-varied moves: (a) $C$ and $C'$ are gradually varied,
(b) $C$ and $C1$ are gradually varied with a difference of $A$ in a 2-cell, (c) $C1$ and $C2$ differ by only two edges in a 2-cell, and
(d)  $C2$ and $C'$ differ only by some edges in a 2-cell. }
\end{center}
\end{figure}

\noindent Therefore, the gradually-varied move from $C$ to $C'$ can be replaced by three simple moves: $C$ to $C1$,  $C1$ to $C2$,
and   $C2$ to $C'$. That is to say that a gradually-varied move from $C$ to $C'$  can be made by a sequence of simple gradually-varied moves
where each adjacent pair only differs by a 2-cell (meaning that only some edges are different  in this 2-cell). This simple gradually-varied move
is called {\it minimal gradual variation} between two curves.

(Both gradually varied moves and minimal gradually-varied moves can naturally be extended to $k$-manifolds. Now we can say that two $k$-manifolds $M$ and $M'$ are gradually varied
if and only if ${\rm XorSum}(M,M')$ is the union of several $(k+1)$-cells. This means that ${\rm XorSum}(M,M')$ contains all 0-cells of these $(k+1)$-cells. We also assume that
$M$ and $M'$ are simply connected. )

In fact, for
any two curves $C_0$ and $C_1$ that are gradually varied, there must be a sequence of curves
where each adjacent pair of curves has the property of minimal
gradual variation.  This means that we only change one 2-cell between two curves
at a time. We can reach $C_1$ from $C_0$ using a sequence of curves.

\begin{lemma}
A gradually-varied move between two discrete curves is equivalent to
 a sequence of discrete curves where
two adjacent curves only change one 2-cell between them
(in terms of ${\rm XorSum}$). This property is true for $k$-manifolds.
\end{lemma}

\begin{proof}
Using mathematical induction to prove the necessary condition, the sufficient
condition is satisfied naturally. First, if ${\rm XorSum}(C_0,C_1)$ only contains one 2-cell,
the condition is valid. Second,
we assume that there are $i$ $2$-cells in ${\rm XorSum}(C_0,C_1)$, which is equivalent
to a sequence of gradually-varied moves where each move only changes one 2-cell.  Then, we
want to prove that if there are $(i+1)$ $2$-cells in ${\rm XorSum}(C_0,C_1)$, then we can still split the moves into
a sequence of gradually varied moves where each move only changes one 2-cell.

Let there be $(i+1)$ $2$-cells in ${\rm XorSum}(C_0,C_1)$.  Since there are a finite number of 2-cells in ${\rm XorSum}(C_0,C_1)$,
  we can select the last 2-cell, denoted by $A$, in ${\rm XorSum}(C_0,C_1)$. We can then construct $C'_1$
such that the boundary of $A$ is ${\rm XorSum}(C'_1,C_1)$.
The rest of $C'_1$ is the same as $C_1$.
We can do this because there are a finite number of cases, and each cell has a finite number of points and edges in discrete space.  Therefore,
${\rm XorSum}(C_0,C'_1)$ contains  $i$ $2$-cells and,
according to the assumption of mathematical induction, we have proven this lemma.
\end{proof}

Again, the lemma states that we can move one curve to another curve gradually by changing one 2-cell at a time.
If we are dealing with a $k$-manifold, then a $(k+1)$-cell is considered. The ${\rm XorSum}$ of the two adjacent
$k$-manifolds (for such a sequence in the proof of the lemma) contains only the boundary of one $(k+1)$-cell.
We provide related examples in the next subsection.




\subsection{Gradually Varied Deformation in a Single Discrete $k$-cell}

In this subsection, we discuss gradually varied deformation of curves and manifolds in a single discrete $k$-cell.
We also show specific examples and methods that can move a discrete curve (or a submanifold) in a $k$-cell.
We are especially interested in moving a $(k-2)$-manifold on the boundary of a discrete $k$-cell, called a $(k-1)$-cycle. This move
is a minimal gradually-varied move. We also note that,
as we discussed in Section 2, the boundary of a $k$-cell is always a $(k-1)$-cycle that is homeomorphic to a $(k-1)$-sphere.
Such homeomorphism is also constructive.


We first look at two examples. In these two examples, we move a curve along the partial boundary of a 3-cell (passing
one 2-cell at a time). This partial boundary of a 3-cell was created by removing a 2-cell from the boundary of the 3-cell.
We need to use this special construction in the following section.

{\bf Example 1: Simplices and Gradually Varied Deformation}

The following example shows how we generate a sequence of curves surrounding the boundary surface of a 3-cell. This sequence
shows the gradual variation of curves. The two ends (curves) of the sequence were originally two closer curves that
only differed by one 2-cell (having the property of minimal gradual variation).

\begin{figure}[h]

\begin{center}

   \epsfxsize=5in
   \epsfbox{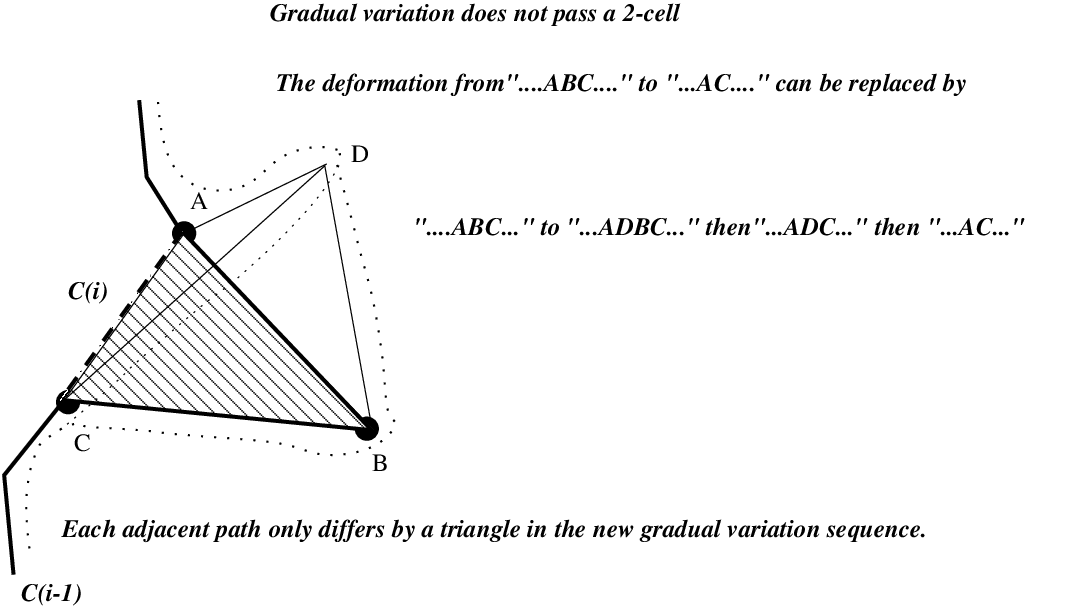}

\caption{Deform one curve to another with gradual variation on the boundary of a 3-simplex. This deformation does not pass the shaded 2-cell. \label{marker}  }
\end{center}

\end{figure}

See the shaded area in Fig.\ 10.
We request that the sequence not pass the shaded area. The question
is as follows: can we always find such a sequence that starts at curve $\cdots ABC\cdots$ and ends at $\cdots AC\cdots$?
In this example, we can first select  $\cdots ADBC\cdots$ and then select  $\cdots ADC\cdots$ to arrive at  $\cdots AC\cdots$.
In other words, without passing the shaded cell, we can use the sequence of $\cdots ABC\cdots$, $\cdots ADBC\cdots$, and  $\cdots ADC\cdots$
to  finally get   $\cdots AC\cdots$.

This is because if we remove a 2-cell from the boundary surface of a 3-cell, the rest of the boundary surface still consists of 1-connected 2-manifolds.
The boundary of the new surface is a 1-cycle.

{\bf Example 2: Cubical Cells and Gradually Varied Deformation}

We now present another example. For a 3-cube, we can still find a sequence that is gradually varied without passing
the bottom 2-cell $ABCD$. (See Fig.\ 11.)

In Fig.\ 11, we have two curves $C_{i-1}=\cdots ABC\cdots$ and  $C_{i}=\cdots ADC\cdots$. $C_{i-1}$ and $C_i$ are gradually varied
with ${\rm XorSum}(C_{i-1},C_i)=$ 2-cell $ABCD$. If the 2-cell is not in consideration, meaning that we cannot
pass this cell for gradual variation, then what we can do is go around the boundary surface of the 3-cell.

The thin-dashed curve in Fig.\ 11 (a) is the first one leaving $C_{i-1}$, which is gradually varied, to go to $C_{i-1}$.
Then, the thin-dashed curve of Fig.\ 11 (b) is the next curve to do the same. If we continue, we will eventually reach $C_i$.

Another property of the sequence is that any two adjacent curves differ by a 2-cell (i.e.,
the only difference between any two adjacent curves is that they share
a 2-cell). There are five steps to get to $C_i$ by passing all five 2-cells.

In higher dimensional cases, we can always make such sequences for simplices or cubic cells.

The easiest case is the $k$-simplex. Let the $(k-2)$-cell be a simple path, where $C_i$ and $C_{i-1}$ share a $k-1$-cell $D$ (${\rm XorSum}$). There
are $(k+1)$ $(k-1)$-cells in the simplex. Removing $D$, we have a $(k-2)$ dimensional boundary set. There is a 0-cell that will be combined with each
of the cells to form a $(k-1)$-cell. Then, we have a total of $k$ cells and a connected $(k-1)$-manifold.

\begin{figure}[h]
\begin{center}

   \epsfxsize=3.2in
   \epsfbox{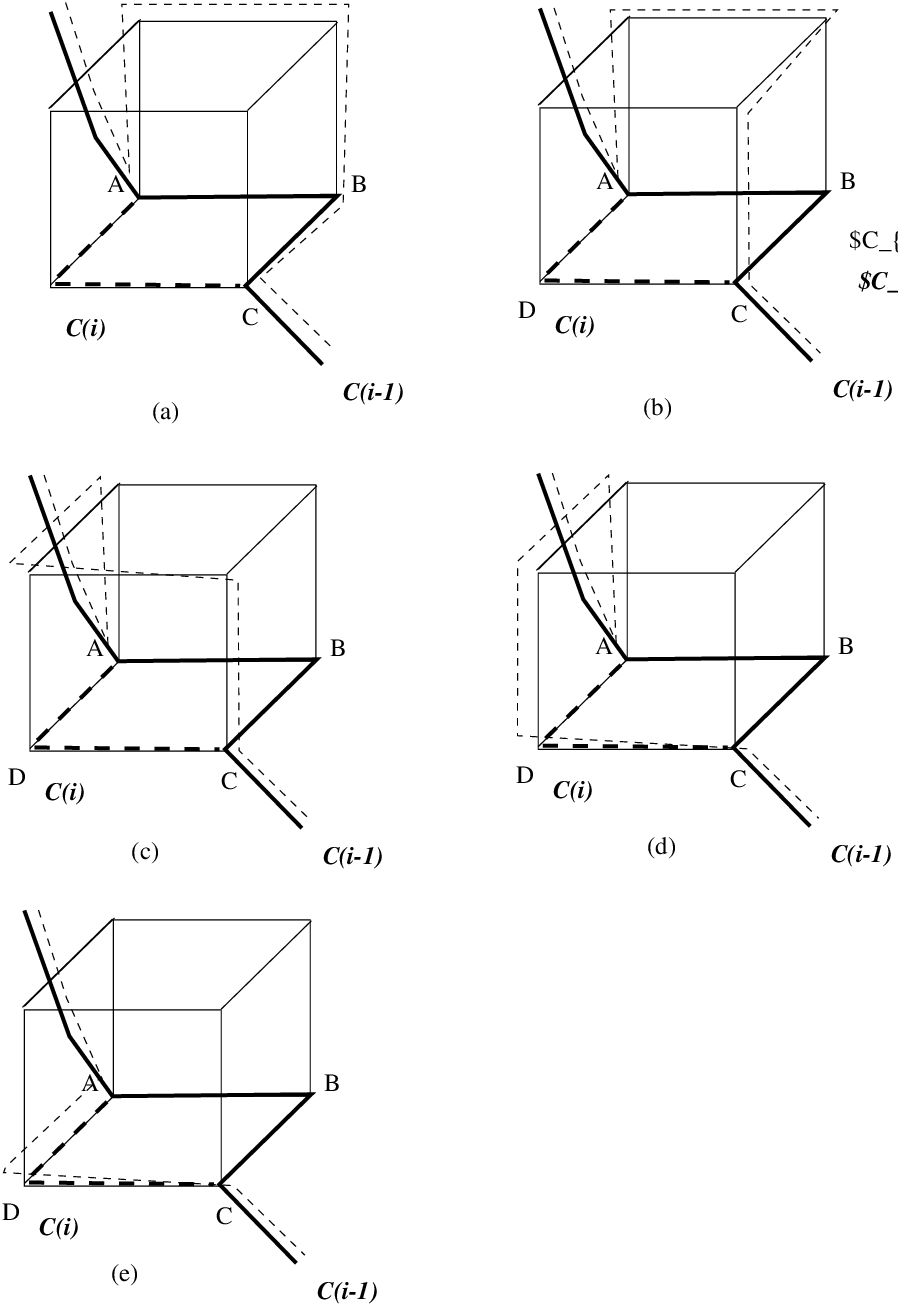}

\caption{Deform one curve to another gradually on the boundary of a 3-cube without passing the bottom 2-cell.  }

\end{center}
\end{figure}

Generally, we can always make such an arrangement.
The boundary of the $k$-cell is a $(k-1)$-cycle, which is (constructively) homeomorphic to the $(k-1)$-sphere.
Removing a $(k-1)$-cell from this boundary will result in creating half a sphere. This half sphere has a $(k-2)$-cycle as its boundary that is homeomorphic to the $(k-2)$-sphere.
In Section 2, we defined $i$-cells for any integer $i$ to be constructive: A cell is defined in the simplest way in terms of being homeomorphic to a ball.
The simplest method is constructive or algorithmic, meaning that we can use an algorithm to determine the path (as well as its boundary, which is homeomorphic to a sphere).

Therefore, we have a generalized conclusion as follows:

\begin{proposition}
    Let $C_0$ and $C_1$ be two discrete $i$-manifolds (curves) with boundaries that are two $(i-1)$-cycles. Let the following three conditions hold:
     (1) $C_0$ and  $C_1$ are joined only by their two end points or at $(i-1)$-cycles,
     (2) The union of $C_0$ and $C_1$ is an $i$-cycle that is homeomorphic to an $i$-sphere, and
     (3) This $i$-cycle is the boundary of an $(i+1)$-cell, $E_{i+1}$, and there is an $(i+2)$-cell $F_{i+2}$ containing $E_{i+1}$ with the boundary that is an $(i+1)$-cycle homeomorphic to an $(i+1)$-sphere.
     Then, we can always find a minimal gradually varied sequence from $C_0$ to $C_1$  without passing $E_{i+1}$.
\end{proposition}

\begin{proof}

We can assume $k=i+2$.  Now, $\partial(F_{i+2}) -E_{i+1}$ means to remove a $(k-1)$-cell, and $\partial(F_{i+2}) -E_{i+1}$ has a boundary that is
a $(k-2)$-cycle. This $(k-2)$-cycle is the union of $C_0$ and $C_1$. We can get a minimal gradually varied sequence from $C_0$ to $C_1$ on $\partial(F_{i+2}) -E_{i+1}$ .

We can use the following method to get such a minimal gradually varied sequence from $C_0$ to $C_1$ on $\partial(F_{i+2}) -E_{i+1}$.
Let $H=\partial(F_{i+2}) -E_{i+1}$.   Cut one $(i+1)$-cell from $H$ at a time (from the boundary of $H$). (This is the inverse process to the connected-sum of adding one $(i+1)$-cell to an $(i+1)$-manifold
with common $(i)$-cell(s). )  In other words, find an $i$-cell $e$ in $C_0$ that is contained in an $(i+1)$-cell $A$ in $H$. $\partial (A) = (\partial (A) -C_0) \cup  ((C_0- \partial (A))\cap \partial (A))$
is an $i$-cycle. Using $(\partial (A) -C_0)$ to replace  $(C_0-\partial (A))\cap \partial (A)$ in $C_0$ will create another $i$-manifold $C_0^{(1)}$. Here,  $C_0^{(1)}=  (\partial (A) -C_0) \cup  ((C_0- \partial (A))$.
We mark the cell as $A$. $C_0$, and $C_0^{(1)}$ are gradually varied with the only difference of passing $A$. Starting at the new $C_0^{(1)}$, we can find another $(i+1)$-cell $A^{(1)}$. Then, we would get $C_0^{(2)}$ and mark a new $A^{(1)}$.
When all $(i+1)$-cells in $H$ are marked, we would have a sequence of gradual variation from $C_0$ to $C_1$. The marking process is to eliminate one cell from $H$ at a time. Since we only have a finite number of cells in $H$, this process will eventually stop.
\end{proof}

%



\section{Discrete Proof of the General\newline Jordan-Schoenflies Theorem}


%
In order to construct a discrete proof of the Jordan-Schoenflies theorem,
we recall some concepts from the theory of discrete manifolds
presented in Section 2. A discrete manifold is a
piecewise linear manifold. The difference is that we cannot
arbitrarily decompose a discrete cell into pieces. For
instance, a discrete 2-cell is predefined. It is a simple and
minimal cycle of discrete 1-cells. Inductively, the boundary of a discrete
$k$-cell is a simple and minimal cycle of discrete
$(k-1)$-cells. Therefore, in our discrete geometry, a discrete
manifold is defined on a graph with topological structures. In
addition, it is finite.

Essentially, in this paper, we do not allow a cell to be
decomposed into smaller cells unless it is expressed
explicitly. In such a discrete case, a cell is already the
minimal entity in its dimension. A $k$-cell is a
$k$-polyhedron but cannot be decomposed into smaller pieces.

The general Jordan-Schoenflies theorem states that, in discrete space, every closed and simply
connected $(n-1)$-submanifold $S$ with local flatness in a closed and simply connected $n$-manifold decomposes
the space into two components, and $S$ is their common boundary. Each of the two components is homeomorphic to the discrete $n$-ball.
In other words, an $(n-1)$-sphere can be locally flatly embedded into an $n$-sphere as its equator.

In \cite{Che13}, we used these discrete techniques to prove the classical Jordan curve theorem:
A closed, discrete curve $C$ separates the plane into two components.

If $M$ is a closed, two-dimensional surface, then $M \setminus C$ consists of two connected components.
We also proved that: If we select a point not on $C$, then there is a component that contains a finite
number of 2-cells, and this component  (when embedded in Euclidean space) is homeomorphic to
a disk.  For a closed $M$ in discrete form (or a piecewise linear 2-complex),
both of the components contain finite numbers of 2-cells (determined by minimal cycles).
Then each of them is homeomorphic to a disk using the same proof.


Our proof in this paper is based on the original, classical Jordan-Schoenflies theorem. In other words, we admit the Jordan-Schoenflies theorem for simply connected closed discrete 2-manifolds (or piecewise linear 2-manifold) $M$:
A 1-cycle that is a discrete curve (which is not a minimal cycle) divides $M$ into two components. Each component is homeomorphic to a 2-cell.

We would also like to restate that if a 1-cycle is a minimal cycle, then this
cycle might be the boundary for a 2-cell in discrete space. A 2-cell in discrete space cannot
be divided into other 2-cells based on our definition (in discrete space). We reject such
a case in order to preserve the properties of the original Jordan
curve theorem. In addition, the union of two 2-cells in this paper is not considered a 2-cell in this paper.  All
2-cells are pre-defined in the discrete case,
but the union of two 2-cells with a common edge will be homeomorphic to
a 2-cell in Euclidean space.   We also assume that $M$ is orientable.

Our proof is divided into two parts: (1) We prove the
Jordan theorem for a closed surface on the 3D manifold, and then
(2) we prove the Jordan-Schoenflies Theorem on the 3D manifold.

\begin{theorem}[Jordan Theorem for a closed surface on the 3D manifold]   \sl
Let $M$ be a simply connected 3D manifold
(discrete or piecewise linear); a closed discrete surface $S$ (with local flatness) will separate $M$ into two components. Here, $M$ can be closed.
\end{theorem}

\noindent {\bf Proof:}  For the beginning part of the proof,  we use the idea of the  proof of the classical Jordan Curve Theorem in \cite{Che13}.
However, the proof given here is independent.


Select a 2-cell in $S$.
This 2-cell is contained in two 3-cells in $M$, called $A$ and $B$. Let $a$ and $b$ be two points in $A$ and $B$, respectively.
Both $a$ and $b$ are adjacent to the intersection of $A$ and $B$.  This intersection is the original 2-cell we chose in $S$.

We know that $a, b \in M \setminus S$. We also know there is
a path $P(b,a)$ from $b$ to $a$ ($P(b,a)$))
passing through a point in $A\cap B$. We denote this by $P^{-1}(b,a)$, the  reverse order of $P(b,a)$.

Now, we want to prove that every path $P(a,b)$ from $a$ to $b$ includes a point in $S$.

On the contrary, we assume there is a $P(a,b)$ that does not include any point in $S$.
Since $M$ is simply connected, there will be a sequence of simple paths (pseudo-curves), $P(a,b)=P(0), P(1),\cdots, P(n)=P^{-1}(b,a)$,
that are side-gradually varied to $P(b,a)$. Note that  $P(b,a)$ contains a point in $S$. There must be a first $i$ such that $P(i-1)$ does not contain
any point in $S$, but $P(i)$ contains $x \in S$. Note that $P(i-1)$ and $P(i)$ are side-gradually
varied, meaning that ${\rm XORSum}(P(i-1),P(i))$ are a collection of 2-cells in $M$.  ${\rm XORSum}(P(i-1),P(i))$ is the exclusive sum that
contains 2-cells where all corner points in each 2-cell are contained in the edges in $P(i-1)$ or $P(i)$ but not in both.
(See the Appendix.)

We illustrate this in Fig.\ 12 below for the current proof, $k=3$:

\begin{figure}[h]

	\begin{center}

   \epsfxsize=5in
   \epsfbox{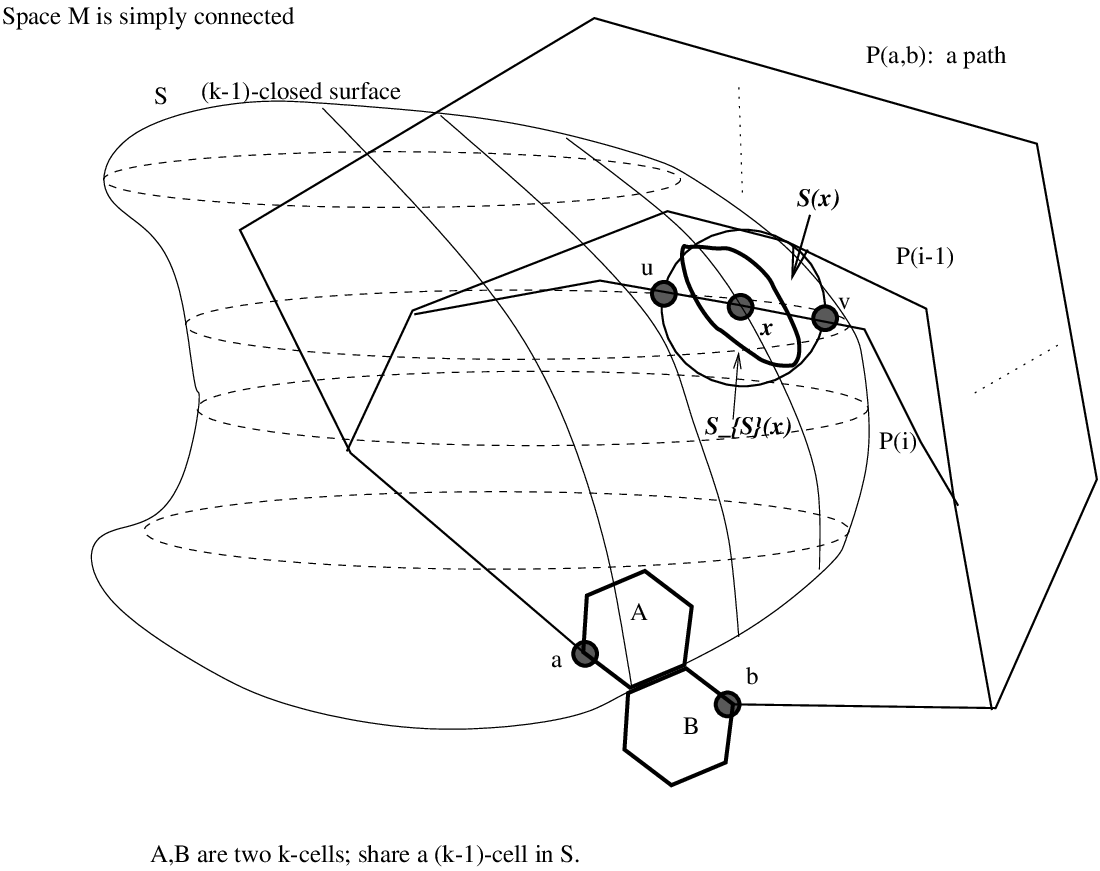}   

\caption{A closed $(k-1)$-manifold in a $k$-manifold. Assume that $x$ is the first intersection point in $S$
where $x\in P(i)$, the collection of paths $P(i)$, $i=1,\cdots,n$. \label{fig_old_1}}

\end{center}
\end{figure}

For a vertex $x$, we use $S_{M}(x)$ to denote a set that contains all 3-cells in $M$ containing $x$, and $S_{S}(x)$
is a set that contains all 2-cells in $S$ containing $x$.
We know that $S_{S}(x) = S\cap S_{M}(x)$.

$S_{M}(x) \setminus \{x\}$ is a 2-sphere, a special 2-cycle where each point is included in a 3-cell that contains $x$.
It is similar to the set where each element links to the center $x$ by one edge or is in the same 3-cell containing $x$.

Note that for higher-dimensional cases, we defined the $k$-cell-distance as the length of the shortest path of $k$-cells (see Section 2).
Now, $k$-cell-distance 1 means that each element $x'$ in the $(k-1)$-cycle,
$S_{M}(x) \setminus \{x\}$,
is in a $k$-cell that also contains the center $x$ ($x'$ is not the center $x$).

In this proof, $S_{S}(x) \setminus \{x\}$ is a 1-sphere. Also, $S_{S}(x) \setminus \{x\}$ is a subset of $S_{M}(x) \setminus \{x\}$.
In other words, $C=S_{S}(x) \setminus \{x\}$ is a closed discrete curve of $Q=S_{M}(x) \setminus \{x\}$.

In the path $P(i)$, there must be a node $u$ that moves to $x$ and also a node $v$ that comes after $x$. See Fig.\ 12.
We have
two major cases: {\bf (a)} $v$ is not on $S$, and {\bf (b)} $v$ is on $S$. We can see that if $v$ is on $S$, then we need to apply
the property of local flatness (discussed in Section 3) .\\
\smallskip \\

\noindent {\bf Proof of Case (a):}
\smallskip \\

According to the Jordan curve theorem,  $S_{S}(x) \setminus \{x\}$ divides $S_{M}(x) \setminus \{x\}$ into two components.
In Fig.\ 12, $\cdots u\to x\to v\cdots$ is a substring (path) of $P(i)$ . Note that $u$ and $v$ are in $Q=S_{M}(x) \setminus \{x\}$ with all
vertices surrounding $x$ in $M$.
Also, $u$ and $v$ are not in $C$ because $C \subset S$. (This means that the path $\cdots u\to x\to v\cdots$ is
 transversal to $S$. We later discuss the case where $v$ is in $S$ so we can find $v'\in Q$.)

Our purpose is to show that from $u$ to $v$, there is a path that shares a part of $P(i-1)$. And this part
must contain a point in $C\subset S$. This generates a contradiction that $P(i-1)$ does not contain any point
in $S$.

Since $P(i-1)$ and $P(i)$ are gradually varied, $u$ will be included in a 2-cell containing
a point in $P(i-1)$, and $v$ will be included in a 2-cell containing a point
in  $P(i-1)$.  Particularly, $x$ will be included in a 2-cell that contains a point in $P(i-1)$.
There will be a cycle $u,\cdots, y_1,\cdots,y_t,\cdots, v,x,u$ that contains part of the path $P(i-1)$,
denoted by $y_1,\cdots,y_t$. This cycle is in $S_{M}(x)$ (Also see an example in \cite{Che13}).
All points $u,\cdots,u_0, y_1,\cdots,y_t, v_0\cdots, v$ are in $Q$.
In fact, this path must contain a point in $C$ based on the Jordan curve theorem. Otherwise, $P(i)$ is not cross-over (or transversal to) $S$.

This path has three parts: $u,\cdots,u_0$ are in $P(i)$; $y_1,\cdots,y_t$ are in $P(i-1)$; and $v_0\cdots, v$ are in $P(i)$. We want to
prove that only the second part $y_1,\cdots,y_t$ can intersect $C$.

We now just need to check whether $u,\cdots,u_0$ could contain any point in $C$.  Since $u_0,\cdots,u$ is on $P(i)$ and
$u_0,\cdots,u$ are in a 2-cell, we see that $x$ is the first element of $P(i)$  on $S$. This is impossible since $C$ is a subset of $S$.

Next, we check whether $v_0,\cdots, v$ could contain a point in $C$.  Note that $v_0,\cdots, v$ are also
in $P(i)$ and they are in a 2-cell including $x$. If $v_0\cdots, v$ has one point in $C$, then this 2-cell contains two or more points in $S$. This is
impossible for the following reason:  Let us say that $v_i$ is in $C$, the path joining with $C$ is
$x\in S, v, \dots, v_i\in C, \dots v_0$. We know that $(x,v,\dots,v_0)$ is contained in a 2-cell denoted by $A1$.
Since $v_i$ is also in $C$, there must be a cell $B1$ in $S_{S}(x)$ that contains both $x$ and $v_i$. ($C=S_{S}(x) \setminus \{x\}$ is the link of $x$.)
The set $A1\cap B1$ contains $x$ and $v_i$ but not $v$.
In the Appendix, we have strictly specified that
any intersection must be a connected path.

(In our definition of discrete manifolds, any two cells must be well-attached or not attached. In other words,
the intersection must be a simply connected  $i$-manifold composed of $i$-cells and be homeomorphic to an $i$-ball when embedding to ${\mathbb R}^n$.
In terms of cell-complexes, the intersection is an $i$-cell that is homeomorphic to an $i$-ball. However, it is hard to determine this fact in
continuous space, and we only mention it because it is not computable.
 It is obvious that
we do not want to allow a complex case of the intersection of two cells. For us,
any two cells $A$ and $B$ can be in any dimension.  The intersection of $A$ and $B$ is a simply connected $i$-manifold and the intersection is homeomorphic to the $i$-ball.
It is important that these facts can be determined in polynomial time of $O(|A|+|B|)$. Here, $|A|$ refers to the number of vertices in $A$.)

In other words, the intersection of two 2-cells or any two cells must be connected by its vertices. However, $v$, \dots, $v_{i-1}$ are not in the intersection. So,
$v_i$ must be $v$. We already assumed that $v$ is not on $S$.
Therefore, there must be a $y_i$ in $y_1,\cdots,y_t$ that is in $C$. Thus, we have a contradiction.
We have proven the case of $S\cap P(i) = \{x\}$ where $x$ is a simple point. \\
\smallskip \\


\noindent {\bf Proof of Case (b):}
\smallskip \\

If $S$ contains two points of $P(i)$, then $C =S_{S}(e=(x,x_1))$ and $G$ are still cycles.
(Any {\it link} of a $k$-face or $k$-cell is a cycle or sphere by a standard theorem in intersection homology theory\cite{GM}). We can still prove the same result as in Case (a).

({\bf Note:}  We require that $S$ be a discrete 2-manifold or $(k-1)$-manifold in $M$. There will not be a case
where $\cdots x v \cdots$ in $P(i)$ if $x$ and $v$ are in $S$, but the edge $(x,v)$ is not in $S$ when edge $(x,v)\in M$. We only allow pseudo-manifolds while we are doing a contraction or other actions; we do not allow a pseudo-manifold when we first select it. Please see the Appendix for more details.
This restriction is related to the so called partial graph properties, meaning that if the vertex set is determined, then the subgraph will
contain all edges (1-cells and $i$-cells) if these vertices are in the set. So, the definition of discrete manifolds here will give
a unique interpretation.    We know that this
situation could also be prevented by the local flatness of $S$ in $M$. We now return to our proof.)

Here is a complex example of this case: $S$ contains a consecutive part of $P(i)$, $X=\{x,x_k,\dots,x_0\}$, i.e., there are more
than two points in $P(i)$ that are in $S$.
Please note that $P(i)$ is a simple path, so $X$ is not restricted to be a flat path. The key of the proof is to
modify $P(i)$ to be a flat path.

Let us have more detailed explanations. In $2D$, $S$ is a curve. So, $x,x_k, \dots x_0$
is a subset of $S$ that is a pure discrete curve. In this paper, $S$ is a surface. $P(i)$ can be zigzaged on $S$.
we may have some cases where $S_{S}(X) \setminus X$ is not a simple (discrete) closed curve.
 Here, $x,x_k,\dots,x_0$ is only a subset of a simple path in higher-dimensional space.

In order to treat this case, we want $x,x_k,\dots,x_0$ to have a collar, meaning that the neighborhood of $x,x_k,\dots,x_0$
does not intersect itself, which is the concept of local flatness. (See Proposition 1.)
In other words,   $S_{S}(X) \setminus X$ and $S_{M}(X) \setminus X=S(X) \setminus X$ must be
a (simple) 1-cycle and 2-cycle, respectively.

We know that $S$ is locally flat in $M$ by the condition of the generalized Jordan-Schoenflies theorem.
This means that $S$ is not folding together in $M$ (see examples in Fig.\ 3 and Fig. 4.). If it does, then we can never make a
locally flat $P(i)$ in $M$.  Therefore, this is our pre-condition.

The key to the proof of Case (b)
is to design an algorithm that uses a technique to modify $P(i)$ into a local flat path:
If $X$ contains two points $a$, $b$ (these two points are not adjacent in $P(i)$) in $S$ such that $a$ and $b$
are adjacent in $S$ or $d(a,b)=2$
(see Fig. 13), then we can have a $P'(i)$ that is gradually varied to $P(i)$. $P'(i)$ still
contains $x$
as the entering point.

Before we describe this modification algorithm, we recall a little more about local flatness.
To observe that $P'(i)$ is ``locally flat'' in $S$  means that there are at least a distance of 3
between any two points that
are not adjacent in the path
$X$ or consecutive such as $a, b, c$ (see Definition 3.1 and Definition 3.2). 

\begin{figure}[h]
	\begin{center}

   \epsfxsize=5in
   \epsfbox{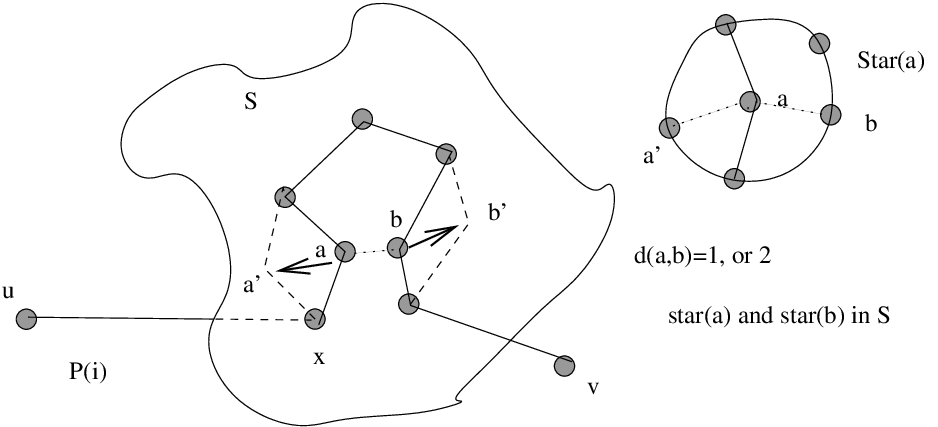}

\caption{ $X \in P(i)$ is not locally flat. We can modify $X$ to make a new locally flat $P'(i)$ in $S$. \label{fig_old_2} }

\end{center}
\end{figure}

Now we describe the algorithm for  modification as follows:

The idea is to insert a sequence of side-gradually varied paths between $P(i-1)$ and $P(i)$. This sequence does not
contain $a$.
In more detail,  we want to find a locally flat $P'(i)$ containing the original $x$ just before reaching $P(i)$ (from $P(i-1)$).

Here, the path $P'(i)$ is gradually varied (by a sequence of paths) to  $P(i-1)$.  In other words, we have a sequence of side-gradually varied paths
from $P_{i-1}$ to $P'(i)$. Except for $P'(i)$, any path in the sequence will not contain a point in $S$. In addition,
this sequence does not contain point $a$. Let us prove this statement: First we draw Fig. ~\ref{fig_old_3},
which is
a continuation of the case of Fig. ~\ref{fig_old_2}.

\begin{figure}[h]
	\begin{center}

   \epsfxsize=5in
   \epsfbox{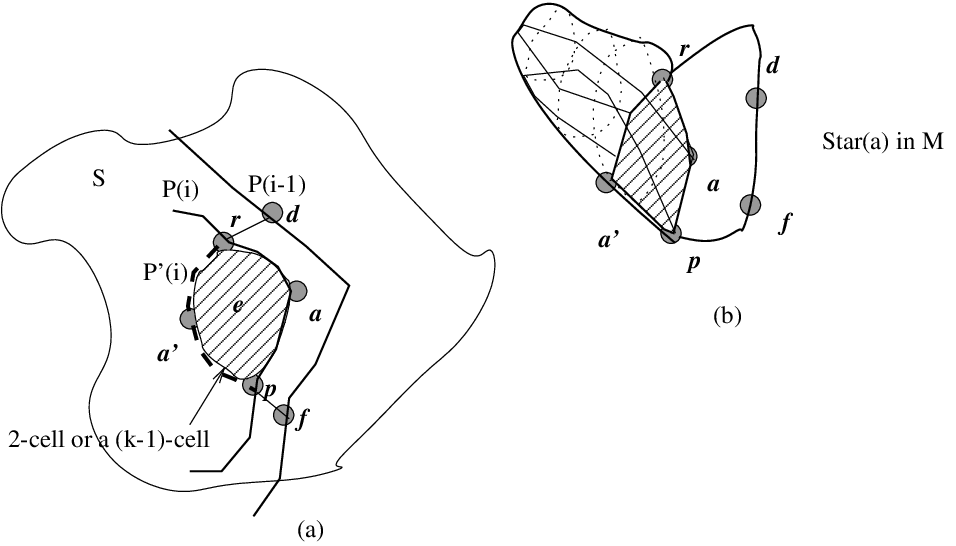}

\caption{Some facts about $P(i)$ and  $P'(i)$: The 2-cell $e$ is between $P(i)$ and  $P'(i)$.
When $a$ moves to $a'$, path $P(i)$ changes to path $P'(i)$, and arc $rap$ changes to arc $ra'p$.
Cell $e$ containing $rapa'$ is in $S$, but $d$ and $f$ are not; $d$ and $f$ are in $Star(a)$. There are two cases:
(i) A 3-cell (or $k$-cell) contains $e$ and arc $df$, and (ii) A 3-cell contains $e$ and is 2-connected to the cell containing arc $df$.
Our purpose is to make a gradually varied path-sequence on $\hbox{\rm Link}(a)$ from $P(i-1)$ to $P'(i)$ where
every path in the sequence (except $P'(i)$) does not intersect $S$. \label{fig_old_3} (Without pass $e$, see Proposition 4.2). \label{fig_old_3}}

\end{center}
\end{figure}

The actual procedure to find such a sequence of side-gradually varied paths is the following:

We want a locally flat $P'(i)$ to replace  $P(i)$ (this process may take more than one iteration).  The point $P'(i)$  is obtained by moving the point $a$ by a 2-cell to
point $a'$ in $S$ to build distance from point $b$ in Fig. ~\ref{fig_old_2}. (If $S$ is $(k-1)$-dimensional, we will move $a$ along with a $(k-1)$-cell.)

In Fig.\ 14,
we exhibit the relationship between $P_{i}$ and $P'_{i}$. The 2-cell $e$ is between these two paths, and
we want to find another path  from $P_{i-1}$ to $P'_{i}$ without passing through $e$.
In other words, we want to find a sequence of side-gradually varied paths from $P_{i-1}$ to $P'_{i}$ that do not pass the 2-cell $e$ (
in Fig. 14).
This is definitely possible since $e$ is contained in a 3-cell, and the inserted paths can go by way of other edges (or faces) to reach  $P'_{i}$ from $P_{i-1}$.

(Note that in the current 3D manifold $M$, this distance is
the edge-distance. If $M$ is a $k$-manifold, then $S$ will be a $(k-1)$-manifold.  The point $a$ will move to be the point $a'$
in the next $(k-1)$-cell, at least to be shared with the next $(k-1)$-cell. We can define $d^{i}(x,y)$ as the
$i$-cell distance
for $x$ and $y$ and $x\neq y$, as the smallest number of $i$-cells in the $(i-1)$-connected path
where $x$ and $y$ are at the end of path. The graph-distance means the distance of edges in the path.
In this paper, the edge-distance, or 1-cell distance, is the graph-distance.)

Here are the facts: (1) $P(i-1)$ is side-gradually varied to $P(i)$.  There must be points $r$ and $p$
on $P(i)$ that have adjacent points $d$ and $f$ in $P(i-1)$ in $M$, respectively. Or they are
in the same 2-cell in $M$. (2) $P'(i)$ is almost the same as $P(i)$ except at point $r$. The
path changes from $r$ to $a'$ and from $a'$ to $p$. (3) We want to build a sequence of paths from $P(i-1)$ to
$P'(i)$ without passing $a$. (That is, do not use the cell $e$. For instance, this is always possible on the boundary of a 3-cell.
By cutting $e$ out, we can still have a bounded surface.) See Section 4.2 and Proposition 4.2.

We also know the following facts: (1) ${\rm Star}_{M}(a)$ contains all points of $r,d,f,p,a'$. (2) Since $S$ is
locally flat, then ${\rm Link}_{M}(a)$ is a 2-sphere that contains all of $r,d,f,p,a'$. In fact, $S$ has a collar in $M$.

Now we want to prove that ${\rm arc}( d\to f)$ is a deformation of ${\rm arc}(d\to r\to a'\to p \to f)$.  This must be true
since we already proved in \cite{Che13}: A discrete 2-cycle (or 2-sphere)\footnote{We can just view
the $k$-sphere as ${\rm Star}(x) \setminus \{x\}$} is simply connected using graph-distance
for contraction. The problem here is finding a path that is side-gradually varied to $P'_{i}$ on ${\rm Link}_{M}(a)$
that does not contain any point in $S\cap {\rm Link}_{M}(a)$.  (This is because there may be other points between $r$ and $a'$.)

Here is the process that can make such a path be the ``collar-edge'' of the ${\rm arc}(r\to a' \to p)$. We know that there
is a side-gradually varied sequence ${\rm arc}(d\to f)=B_1$, $B_2$,\dots,$B_t= {\rm arc}(d\to r\to a' \to p\to f)$. Let $B_i$
be the first containing point $q$ (which is neither $r$ nor $p$) that is in a 2-cell or $(k-1)$-cell containing a point in  ${\rm arc}(r\to a' \to p)$.
We want to fix $q$ in $W$, a queue that was originally empty, so we split $B_i$ into two sub-arcs: $A_1$ from $r$ to $q$, and
$A_2$ from $q$ to $p$. (We will have two smaller 1-cycles based on the two sub-arcs $A_1$ and $A_2$.) Therefore we can repeat this process to make a $q'$
that lies between $r$ and $q$ so that $q'$ shares a 2-cell with a point in $P'_{i}$, and put $q'$ into $W$. Since we have a finite
number of points, we can perform the same process for $A_2$.

So, $dWf$ will be a path that is side-gradually varied to $P'(i)$ in ${\rm Link}_{M}(a)$, where $\rm Link$ is ${\rm Star}(x) \setminus \{x\}$.
Note that $d$ and $f$ in Fig.\ 14
are points on $P(i-1)$, which has an edge-distance of 1 to $P(i)$.
In addition, $dWf$  does not contain any points in $S$.  Replacing $dWf$ in $P(i-1)$, we get $P'(i-1)$.
In other words, we use the arc $rdWfp$ to replace the arc $rd\cdots fp$ in  $P(i-1)$ to get the new path $P'(i-1)$.
$P'(i-1)$ has a gradually varied path sequence to $P(i-1)$ (see Section 4.2 and Proposition 4.2.).

Thus, $P'(i)$  is locally flat at the point we modified. If $P'(i)$ contains cases that are not locally flat, we will need to
repeat the above process by changing another $a$ to $a'$ until the entire $P'(i)$ is locally flat.  This procedure is
finite since our space is finite, meaning that there are only a finite number of cells in our discrete spaces.

There must be a path in the deformation sequence to locally
flat $P'(i)$ (in $S$). (See Section 3 and Section 4.)  The set $X'=\{x, \dots,\}\in P'(i)$ is a
subset of $S$. Because $X'$ is locally flat, we have two
cycles $C^{(1)}=S_{S}(X') \setminus X'$ and $C^{(2)}=S_{M}(X') \setminus X'$, where
$S_{S}(X') \setminus X'$ is a 1-cycle and $S_{M}(X') \setminus X'$ is a 2-cycle
when $M$ is a 3-manifold. $C^{(1)}$ is the closed curve in
$C^{(2)}$. According to the Jordan curve theorem, every path
that is gradually varied to $P'(i)$ must contain a point in
$C^{(1)}$. Assume that $P'(i-1)$ is such a path toward
$P(i-1)$ from $P'(i)$ (we can denote the gradually-varied
paths as $P'(i), P'(i-1), P^{(2)}(i-1),\cdots,
P^{(t)}(i-1)=P(i-1)$), then $P'(i-1)$ must contain a point in
$C^{(1)}$. So, $P'(i-1)$ has a point $x'$ in $S$. Based on
$P'(i-1)$ and $x'$, we can continue the above process until we
reach $P(i-1)$ (this is because we only have a finite
number of paths). By continuously doing this analysis,  $P(i-1)$ must have a point
in $S$. Therefore, we have a contradiction. Thus, we have proven the statement
we wanted for Case {\bf (b)}, which is a key construction.

(For a $k$-manifold $M$, we will have a similar construction.
If $M$ is a $k$-manifold, then $S_{S}(X) \setminus X$ is a $(k-2)$-cycle
and  $S_{M}(X) \setminus X$ is a $(k-1)$-cycle. See Corollary 3.7.)

Thus, we have proven the generalized Jordan-Schoenflies theorem for the discrete closed surface in 3D space.
\endpf
\smallskip \\

We would like to discuss a little more about local flatness even though Section 3 is dedicated to this topic. We can view
 local flatness in the discrete case as follows: A discrete curve $C$ is said to be
{\it locally flat} if for any proper subset (arc) $X$ of $C$, $S(X) \setminus X$ is a simple cycle.
See the proof of Theorem 3.6 where the existence of the collar preserves this property. In addition,
the manifold $M_n$ must be locally flat. Otherwise, $M_k$ cannot be locally flatly embedded in $M_n$.
(We usually view $M_n$ as an $n$-sphere, $n$-cycle, or a manifold that is homeomorphic to an $n$-sphere as predefined and computationally decidable in
polynomial time.)

On the other hand, a locally flat path $P(t)$ means that $S(x)$ ($x$ with a collar) does not intersect
with $S(x')$ if $x,x' \in P(t)$ are apart from each other. Intuitively, the  collar of $P(t)$ is just the union of $S(x)$ for all points $x\in P(t)$.
This means that $x$ to $x'$s in $M$ or $S$ must be apart from each other with a distance of 3.

Therefore, the graph-distance (or cell-distance) of 3 is the key for most pairs of points or cells. See Section 3.2 and Definition 3.2.
Distance 3 is the minimum distance in the discrete case for local flatness where the collar will not intersect,
except in some cases found in Fig.\ 5 in Section 3.
With flatness, ${\rm Link(subpath}={\rm arc}) = {\rm Star(arc)} \setminus \{{\rm arc}\}$ is a cycle in either $S$ or $M$.
\footnote{Graph-distance (or edge-distance) 3 in our previous version posted on Arxiv.org was not very accurate since we did not consider certain
``corner'' cases. In Section 3 of this paper, we expanded further investigations on discrete local flatness. We prove the following statement in the proof of Theorem 3.6 in Section 3:
With local flatness, ${\rm Link(subpath}={\rm arc}) = {\rm Star(arc)} \setminus \{{\rm arc}\}$ is a cycle in either $S$ or $M$. This
is also the real meaning of the collar of a curve.}

Now, we prove the Schoenflies theorem: A closed 2-cycle separates $M$ ($M$ is homeomorphic to 3-sphere) into two components, each of
which will be homeomorphic to a 3-ball.  (We use a technique similar to the one used in \cite{Che13}.)

The following lemma completes the proof
of the Jordan-Schoenflies theorem for a closed surface in the 3D sphere $M$.

\begin{lemma}  \sl
Let $M$ be a simply connected 3-manifold (discrete or piecewise linear) that is homeomorphic to
a 3-sphere.  A closed discrete $(2)$-manifold $S$ on $M$ will separate $M$ into two components where
each component is homeomorphic to a 3-cell.
\end{lemma}


\noindent {\bf Proof:}
Since $S$ separates $M$ into two, we can find a connected component and mark every
3-cell in this component; we denote this component as ${\cal D}$. Choosing one 3-cell $D_p$ that has a 2-face (2-cell in this case) in $S$, we
design a procedure to contract $S$ toward to this 3-cell $D_p$.


The algorithm is similar to the algorithm described in \cite{Che13}. We first measure (compute)
the distance
from each $k$-cell ($k=3$) $X$ in the component ${\cal D}$ to $D_p$. This distance
counts how many $k$-cells are between $X$ and $D_p$ (in other words, the distance is the number of $k$-cells involved). More specifically, note that if $p$ is a point (0-cell) in $D_p$, then, each point $q\ne p$ in $D_p$ has
a distance of 1. Each point $r$, not equal to $q$, in another $k$-cell (not $D_p$)
that contains both $r$ and $q$ will have a distance of 2 to $p$ if $r$ is not in $D_p$.

In other words, use the 3-cell distance to measure how far each 3-cell is to $D_p$ in the component. (See Section 2.4.)
This distance is called the $k$-distance and indicates the length of the shortest path.


As a result, we can find the longest 2-cell $F$ from $S$ in $D$, meaning that there is a point $x\in F$
that has the longest $k$-distance to $p$. This must exist since we only have a finite number of
$k$-cells in the component. Note that $S$ and $M$ are orientable, and so there must also be a longest distance on $S$.
We also assume $F$ is a 2-face in a 3-cell $X$ in ${\cal D}$.

Now, we delete a 2-cell $F$  in $S$ containing $x$.

After we delete this face (2-cell) in $S$, we use other faces in $X$ that are not
in $S$ to replace the deleted one. Thinking about the intersection of $S$ and $X$ ($X$ is a 3-cell containing $F$, $X \cap S$ is not empty. ),
this intersection can be empty, 0-cells, 1-cells, or $(k-1)$-cells (2-cells in our case).
We are only interested in the intersection that is a 2-cell or a set of 2-cells (this set must be homeomorphic to a 2-disk).
If the intersection is the 2-cell $F$, then all 2-cells in $X \setminus F$ will be like a coffee cup without
a lid $F$.

Let $U_2(X)$ be the set of 2-cells in $X$. (We also define $U_k(X)$ to be the set of $k$-cells in $X$.)
$U_2(X) \setminus \{F\}$ will contain all 2-cells in $X$ except $F$, so $(U_2(X) \setminus F) \cap F$ is the boundary cycle of $F$ and $U_2(X) \setminus \{F\}$.

Thus,  $(U_2(X) \setminus F) \cap F$  is a simple closed path (or a closed $(k-2)$-pseudo-curve or manifold).   If $S\cap X$ is a set of 2-cells, then these 2-cells
are connected ($S$ is a pseudo-manifold), and the boundary is also a simple closed
path (or a closed $(k-2)$-pseudo-curve or manifold). Let $B(X)$ be the boundary faces of
$X$. Using $B(X) \setminus (S\cap X)$ to replace $(S\cap X)$ in $S$, we will have an $S'$.

(We do not consider the case where $S\cap X$ does not contain 2-cells or $(k-1)$-cells in this proof.)

We have a new $S'$ that is also a closed
pseudo-surface. ($S$ is a surface.)  The new ${\rm XorSum}(S',S)$
is the boundary of $X$.  Changing $S$ to $S'$, we
will reduce the internal part of $X$, i.e. we determined a 3-cell from the original
component $D$. We can mark it.

In the above process, we removed a 3-cell $X$. We can repeat this process to remove all 3-cells
except $D_p$. This is because the number of (both unmarked and marked) 3-cells is finite, and the process will stop eventually.
This is the general philosophy of the algorithm. When some extreme cases occur, we also need to process it. For instance, when
the 2-cell $F$ found is contained in $X$ where $X$ has a point that is on $S$ but this point is not connected to $F$ in $S\cap X$. 
In other words, $S\cap X$ is not a connected or it is not homeomorphic to a 2-disk. This $X$ could not be deleted by removing $F$ to 
keep ${\cal D}\setminus X$ to be simply connected. In such a case, we just need to remove any $F$ so that $X$ can be removed from ${\cal D}$
to maintain it to be simply connected. This means that $S\cap X$ is homeomorphic to a 2-disk and  $S\cap X$ is not the boundary of a 3-cell. 
In fact, theoretically we only need to
remove those $X$ to keep the new ${\cal D}$ to be simply connected in contracting process.  Such $X$ always exist for ${\cal D}$ 
~\footnote  {An extreme case was found when the first author revising his paper in dealing with 3D triangulation of a 3D compact manifold.
The statement is still true.  
See Appendix A in {\it L. Chen, Algorithms for Deforming and Contracting Simply Connected Discrete Closed Manifolds (III), 2017.}  $https://arxiv.org/pdf/1710.09819.pdf$. (The revision v5 was
posted in Feb. 2020).}. 
%
We can see that $S$ can be contracted to the boundary of a 3-cell; then,
we can contract the 3-cell to a point. The inverse of this process will provide
a homeomorphic mapping from the component bounded by $S$ to a 3-cell, $D_p$.

When we deal with a $k$-manifold, the principle of the proof is the same.
\endpf
\smallskip \\


\begin{theorem}[The general Jordan-Schoenflies Theorem]  \sl
If $M$ is a simply connected $k$-manifold (discrete or piecewise linear)
that is homeomorphic to a $k$-sphere, then
a closed discrete $(k-1)$-manifold with local flatness in $M$ will separate $M$ into two components.
In addition, each component is homeomorphic to a $k$-cell.
\end{theorem}

\noindent {\bf Proof:}
We use mathematical induction to prove the case for a (discrete or piecewise linear) simply connected $k$-manifold $M$ that is homeomorphic to a $k$-sphere;
a closed $(k-1)$-manifold $S$ will divide $M$ into two components. The assumption is that a $(k-1)$-manifold that is homeomorphic to a $(k-1)$-sphere
satisfies the Jordan-Schoenflies theorem. Just like we did in the proof of Theorem 5.1, if there is a path from $a$ to $b$ without passing any point in $S$, then we can denote the path
as $P=P(0)$.  If $x$ is the first point that is on $S$ in the sequence of side-gradually varied paths $P(0), \cdots, P(n)$,
we can assume that this point is in $P(i-1)$. We can use the same strategy we used before
to prove this theorem: {\bf (a)} If $x$ is the only point in $P(i) \cap S$,  then $Q=S_{M}(x) \setminus \{x\}$ is a $(k-1)$-sphere and $C=S\cap Q$ is a $(k-2)$-sphere.
According to the inductive hypothesis, $C$ separates $Q$ into two components.
Every path from $u$ to $v$ in $Q$ must contain a point in $C$, so $P(i-1)$ must contain a point in $C$.
There is a contradiction if we assume $P(i-1)$ does not contain
any point in $S$. {\bf (b)} If $X=P(i) \cap S$ contains more than one point, then we can use local flatness to find a locally flat $P'(i)$ that is gradually varied to $P(i)$.
In addition, if $X'=P'(i) \cap S$, then $Q=S_{M}(X') \setminus \{x\}$ is a $(k-1)$-sphere and $C=S\cap Q$ is a $(k-2)$-sphere. We can still use the
Jordan theorem for general closed manifolds in $(k-1)$-dimensions to prove the current theorem.

We can use the same technique to prove that the connected component is homeomorphic to a $k$-disk in the above lemma.
Therefore, we have proven the general Jordan-Schoenflies theorem.

\endpf
\smallskip \\


The above Theorem 5.3 can be split into two pieces just like Theorem 5.1 and Lemma 5.2.  The treatment was presented in {\it 2018 Lehigh University Geometry and Topology Conference}.
(L. Chen jointly with S. G. Krantz, The Discrete Method for Decomposition of $n$-Spheres and $n$-Manifolds, Lehigh University Geometry and Topology Conference, May 24-27, 2018.) We just repeat the proof above as follows:

\begin{theorem}  (The general Jordan separation Theorem) If $M$ is a simply connected $k$-manifold (discrete or piecewise linear), then a $(k-1)$-closed discrete manifold, orientable and simply connected with local flatness in $M$ will separate $M$ into two components.
\end{theorem}

\noindent {\bf Proof:}
Using mathematical induction, we can prove that for a (discrete or piecewise linear) simply connected $k$-manifold $M$, a $(k-1)$-closed discrete manifold, orientable and simply connected with local flatness in $M$ will separate $M$ into two components.   The inductive hypothesis here is: For a (discrete or piecewise linear) simply connected $(k-1)$-manifold $M'$, a $(k-2)$-closed discrete manifold, orientable and simply connected with local flatness in $M'$ will separate $M'$ into two components.

Just like we did in the proof of Theorem 5.1, if there is a path from $a$ to $b$ without passing any point in $S$, then we can denote the path
as $P=P(0)$.
If $x$ is the first point that is on $S$ in the sequence of side-gradually varied paths $P(0), \cdots, P(n)$,
we can assume that this point is in $P(i-1)$. We can use the same strategy we used before
to prove this theorem: {\bf (a)} If $x$ is the only point in $P(i) \cap S$,  then $Q=S_{M}(x) \setminus \{x\}$ is a $(k-1)$-sphere and $C=S\cap Q$ is a $(k-2)$-sphere.
According to the inductive hypothesis, $C$ separates $Q$ into two components.
Every path from $u$ to $v$ in $Q$ must contain a point in $C$, so $P(i-1)$ must contain a point in $C$.
There is a contradiction if we assume $P(i-1)$ does not contain
any point in $S$. {\bf (b)} If $X=P(i) \cap S$ contains more than one point, then we can use local flatness to find a locally flat $P'(i)$ that is gradually varied to $P(i)$.
In addition, if $X'=P'(i) \cap S$, then $Q=S_{M}(X') \setminus \{x\}$ is a $(k-1)$-sphere and $C=S\cap Q$ is a $(k-2)$-sphere. We can still use the
Jordan theorem for general closed manifolds in $(k-1)$-dimensions to prove the current theorem.

\endpf
\smallskip \\


\begin{theorem}  (The general Jordan-Schoenflies Theorem) If $M$ is a simply connected $k$-manifold (discrete or piecewise linear) that is homeomorphic to a $k$-sphere, then a $(k-1)$-cycle with local flatness in $M$ will separate $M$ into two components. In addition, each component is homeomorphic to a $k$-cell. (we can assume this $k$-sphere is in Euclidean space.)
\end{theorem}

\noindent {\bf Proof:}
For a (discrete or piecewise linear) simply connected $k$-manifold $M$ that is homeomorphic to a $k$-sphere, we can get that there is a partition or decomposition of $M$ and each $k$-cell of the partition is homeomorphic to a Euclidean $k$-cell. So we can get the discretization of the cell to be a discrete $k$-cell.  Again, the process for homeomorphism can be done in finite time or even in polynomial time regarding numbers of cells in the final decomposed $k$-complex.

The first part of this theorem is the same as Theorem 5.4.  We only need to face the second part.
We can use the same technique to prove that the connected component is homeomorphic to a $k$-disk in Lemma 5.3.
Therefore, we have proven the general Jordan-Schoenflies theorem.

\endpf
\smallskip \\

The advantage of using the discrete method for proving the general Jordan-Schoenflies theorem is rendering
the proof as an algorithmic procedure. We can actually program this algorithm
for contraction. This method may also have applications in other geometric problems.


\section{Conclusion}

In this paper, we give a complete proof of the general Jordan-Schoenflies theorem.
When we say that $M$ is a simply connected $k$-manifold (discrete or piecewise linear) that is also homeomorphic to the $k$-sphere,
we mean that there is an efficient constructive method (a polynomial time algorithm in computational geometry)
to decide whether $M$ is homeomorphic to a $k$-sphere.
We prove the main result at the end of Section 5.

The general Jordan Schoenflies theorem states: Embedding an $(n-1)$-sphere $S^{(n-1)}$ local flatly in
 an $n$-sphere $S^{n}$ decomposes
$S^{n}$ into two components. In addition, the embedded  $S^{(n-1)}$ is the common boundary of the two components and
each component is homeomorphic to the $n$-ball.
According to the theorem in Section 2 that states that local flatness implies the existence of a collar, the general Jordan-Schoenflies theorem can also be stated as:
Embed $S^{n-1}\times [-1,1]$ in the $n$-sphere, then each of the closed components bounded by $S^{n-1} \times 0$ in this  $n$-sphere is homeomorphic to the $n$-ball.



\newpage

\section*{Appendix: Basic Concepts of Discrete Manifolds}
\addcontentsline{toc}{section}{Appendix}

\begin{appendix}

In topology, the formal description of the Jordan curve theorem is:
A simple closed curve $J$ in a plane $\Pi$ decomposes $\Pi \setminus J$ into two
components.  In fact, this theorem holds for any simply connected surface.
A plane is a simply connected surface in Euclidean space, but this theorem is not true
for a general continuous surface. For example, a torus fails this result.

What is a simply connected continuous surface?
A connected topological space $T$ is simply connected
if, for any point $p$ in $T$, any simply closed curve containing $p$ can be
contracted to $p$. The contraction is a continuous mapping among
a series of closed continuous curves~\cite{New}.

In order to keep the concepts simple to understand, we first define the
gradual variation between graphs. Then, we define discrete deformation
among discrete pseudo-curves. And finally, we define the contraction of curves
as a type of discrete deformation.

In this section, we assume the discrete surface is both regular and orientable.
A discrete surface is  regular if every neighborhood of each point is homomorphic to a
2D discrete disk (a umbrella shape)~\cite{Che14}.

\begin{definition}   \rm
Let $G$ and $G'$ be two connected graphs. A mapping
$f: G\rightarrow G'$ is gradually varied if, for two vertices
$a, b \in G$ that are adjacent in $G$,
$f(a)$ and $f(b)$ are adjacent in $G'$ or $f(a) = f(a')$.
\end{definition}

Intuitively, ``continuous'' change from a simple path $C$ to another path $C'$ means that
there is no ``jump'' between these two paths.  If $x, y\in S$, then $d(x,y)$ denotes
the distance between  $x$ and $y$.  For instance, $d(x,y)=1$ means that $x$ and $y$ are
adjacent in $S$. It is important to point out that, in a 2-cell (or any other k-cell),
from a point $p$ to another point $q$ in the cell, $p\ne 1$, the distance $d(p,q)$
can be viewed as $1$. In other words, a cell can be viewed as a complete subgraph
on its vertices.

\begin{definition}   \rm
Two simple paths $C=p_0,\dots,p_n$ and
$C'=q_0,\dots,q_m$ are gradually varied
in $S$ if  $d(p_0,q_0)\le 1 $ and  $d(p_n,q_m)\le 1 $,
and  for any non-end point $p$ in $C$:

(1) $p$ is in $C'$, or $p$ is contained by a 2-cell $A$ (in $G(C\cup C')$)
	  such that $A$ has a point in $C'$.

(2) Each non-end-edge in $C$ is contained by a 2-cell $A$ (in $G(C\cup C')$), which
	 has an edge contained by $C'$ but not $C$ if $C'$ is not a single point.

\noindent And vice versa for $C'$.
\end{definition}

For example, $C$ and $C'$ in Fig.\ 15 (a) are gradually varied,
but $C$ and $C'$ in Fig.\ 15 (b) are not gradually varied.
We can see that a 2-cell, which is a simple path, and any two connected parts
in the 2-cell are gradually varied, so we can say that
a 2-cell can be contracted to a point gradually.


Assume $E(C)$ denotes all edges in
path $C$.  Let ${\rm XorSum}(C,C') = (E(C) \setminus E(C'))\cup (E(C') \setminus E(C))$.
${\rm XorSum}$  is called ${\rm sum ( modulo}\ 2)$ in  Newman's book ~\cite{New}.


\begin{figure}[h]
	\begin{center}

   \epsfxsize=4.5in
    \epsfbox{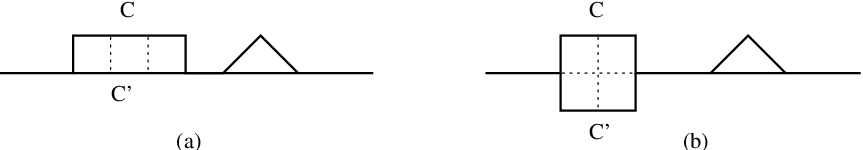} 

(a) \hskip 2in (b)
\caption{Gradually varied curves: (a) $C$ and $C'$ are gradually varied;
(b) $C$ and $C'$ are not gradually varied.}
\end{center}
\end{figure}

Attaching a 2-cell to a simple path $C$, if the intersection is an arc (connected
path) and not a vertex, then we can cut the intersection (keeping the first and last
vertices of the intersection, which is an arc); the simple path will go another
half of the arc of the cell. The new path is also a simple path, and
it is gradually varied to $C$. Therefore,

\begin{lemma}  \sl
 Let $C$ be a  pseudo-curve and $A$ be a 2-cell.
If $A \cap C$ is an arc containing at least one edge,
then ${\rm XorSum}(C,A)$ is a gradual variation of $C$.
\end{lemma}

It is not difficult to see that ${\rm XorSum}({\rm XorSum}(C,A),A)= C$ and
\newline
${\rm XorSum}({\rm XorSum}(C,A),C)= A$ under the condition of the above lemma.

\begin{definition} \rm
Two simple paths (or pseudo-curves) $C$ and $C'$ are said to
be homotopic if
there is a series of simple paths $C_0,\dots,C_n$ such that
$C=C_0$, $C'=C_n$, and  $C_{i}, C_{i+1}$ are gradually varied.
\end{definition}

We say that $C$ can be discretely deformed\index{Discrete deformation}
to $C'$ if $C$ and $C'$ are homotopic. The following
lemma states that we can deform a curve by making changes
one cell at a time.

\begin{lemma} \sl
If two (open, not closed) simple paths $C$ and $C'$ are homotopic, then
there is a series of simple paths $C_0,\dots,C_m$ such that
$C=C_0$, $C'=C_n$, and  ${\rm XorSum}(C_{i}, C_{i+1})$ is a 2-cell excepting
end-edges of $C$ and $C'$.
\end{lemma}


To prove the Jordan curve theorem, we need to describe what the
disconnected components are by distinguishing them
from a simple curve $C$.  It means that any path from one component
to another must include at least one point in $C$.
It also means that this linking path must cross-over the curve $C$.
In this subsection, we want to define this idea.

Because a surface-cell $A$ is a closed path, we can define two
orientations (normals) to $A$: clockwise and counter-clockwise.
Usually the orientation of a 2-cell is not a critical issue.
However, it is necessary for the proof of the Jordan curve theorem.

In other words, a pseudo-curve, which is a set of points with no ``direction,''
as a path has its own ``travel direction'' from $p_0$ to $p_n$.
For two paths $C$ and $C'$,
which are gradually varied, if a 2-cell $A$
is in $G(C\cup C')$, the orientation of $A$ with respect to $C$ is determined by
the first pair of points $(p,q) \in C \cap A$ and $C= \dots p q \dots $.
Moreover, if a 1-cell of $A$ is in $C$, then the orientation of $A$ is
fixed with respect to $C$.

According to Lemma 7.6 in \cite{Che14}, $S(p)$ contains all adjacent points of $p$ and
$S(p) \setminus \{p\}$ is a simple cycle, and there is a cycle containing all points in $S(p) \setminus \{p\}$.

We assume that the cycle $S(p) \setminus \{p\}$ is always oriented
clockwise. For two points $a, b \in S(p) \setminus \{p\}$,
there are two simple cycles containing the path $a\rightarrow p \rightarrow b$:
(1) A cycle from $a$ to $p$ to $b$ then moving clockwise to $a$, and
(2) A cycle from $a$ to $p$ to $b$ then moving counter-clockwise to $a$.
See Fig.\ 16(a).

It is easy to see that the simple cycle $S(p) \setminus \{p\}$ separates
$S \setminus \{S(p) \setminus \{p\}\}$ into at least two connected components because from
$p$ to any other points in $S$, the path must contain a point in $S(p) \setminus \{p\}$.
$S(p) \setminus\{p\}$ is an example of a Jordan curve.

\begin{figure}[h]
	\begin{center}

   \epsfxsize=4.5in
   \epsfbox{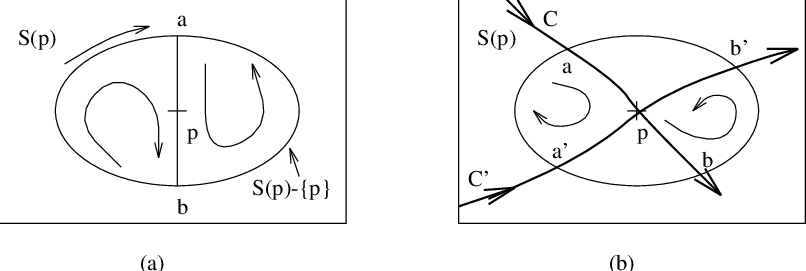}   

\caption{$S(p)$ and cross-over at $p$: (a) Two adjacent points $a$ and $b$ of $p$ in
 $S(p)$, and (b) Two cross-over paths.}

\end{center}
\end{figure}

\begin{definition}   \rm
Two simple paths $C$ and $C'$ are said to  ``cross-over'' each other if
there are points $p$ and $q$ ($p$ may be the same as $q$) such that
$C = \dots a p b\dots s q t\dots$ and $C'= \dots a' p b\dots s q t'\dots$ where
$a\ne a'$ and $t\ne t'$. The cycle
$a p a' \dots a$ without $b$ in $S(p)$ and the cycle $q t \dots t' q$ without
$s$ in $S(q)$ have different orientations with respect to $C$.
\end{definition}

For example, in Fig.\ 16 (b), $C$ and $C'$ do ``cross-over'' each other.
When  $C$ and $C'$ do not ``cross-over'' each other, we will say
that $C$ is on the side of $C'$.

\begin{lemma}	\sl
If two simple paths $C$ and $C'$ do not cross-over each other, and they
are gradually varied, then every surface-cell in $G(C\cup C')$ has the
same orientation with respect to the ``travel direction'' of $C$ and
opposite to the ``travel direction'' of $C'$.
\end{lemma}

We also say that $C$ and $C'$ in the above lemma  are side-gradually
varied\index{Side-gradually-varied}.

Intuitively, a simply connected set is a set where, for any point,
every simple cycle containing
this point can contract to the point.

\begin{definition}  \rm
A simple cycle $C$ can contract to a point $p\in C$
if there exists a sequence of simple cycles, $C=C_0, \dots, p=C_{n}$ such that:
{\bf (1)} $C_i$ contains $p$ for all $i$;
{\bf (2)} If $q$ is not in $C_i$ then $q$ is not in all $C_{j}$, $j>i$;
{\bf (3)} $C_{i}$ and $C_{i+1}$ are side-gradually varied.
\end{definition}
					
We now exhibit three reasonable definitions of simply connected spaces below.
We will provide a proof for the Jordan curve theorem under the third definition
of simply connected spaces. The Jordan theorem shows the relationship
among an object, its boundary, and its outside area.

A general definition of a simply connected space should be:

\begin{definition} {\bf Simply Connected Surface Definition (a)} \rm  $\langle G,U_2\rangle$ is
 simply connected if any two closed simple paths are homotopic.
\end{definition}

If we use this definition, then we may need an extremely long proof for
the Jordan curve theorem. The next one is the standard definition,
which is a special case of the above definition.

\begin{definition}[Simply Connected Surface Definition {\bf (b)}]  \rm
A connected discrete space $\langle G,U_2\rangle$ is  simply connected
if, for any point $p\in S$, every simple cycle containing $p$ can contract to $p$.
\end{definition}

This definition for the simply connected set is based on the original meaning of
simple contraction. In order to
make the task of proving the Jordan theorem
simpler, we give the third strict definition of simply connected surfaces as
follows.

We know that a simple closed path (simple cycle) has at least
three vertices in a simple graph.
This is true for a discrete curve in a simply connected surface $S$.
For simplicity, we call an unclosed path an arc.
Assume that $C$ is a simple cycle with clockwise orientation.
Let two distinct points $p, q\in C$. Let $C(p,q)$ be an arc of $C$ from $p$ to $q$ in
a clockwise direction, and $C(q,p)$ be the arc from $q$ to $p$ also in a clockwise
direction.
Then, we know that $C=  C(p,q)\cup C(q,p)$. We use $C^{a}(p,q)$  to represent
the counter-clockwise arc from $p$ to $q$. Indeed, $C(p,q) = C^{a}(q,p)$.
We always assume that $C$ is in clockwise orientation.

\begin{definition}[Simply Connected Surface Definition {\bf (c)}]  \rm
A connected discrete space $\langle G,U_2 \rangle$ is simply connected
if, for any simple cycle $C$
and  two points $p, q\in C$, there exists a sequence of simple cycle paths
$Q_{0},\dots,Q_{n}$ where $C(p,q)=Q_{0}$ and $C^{a}(p,q)=Q_{n}$ such that
$Q_{i}$ and $Q_{i+1}$ are side-gradually varied for all $i=0,\cdots, n-1.$.
\end{definition}

In fact, it is easy to see that the Definitions {\bf (b)} and {\bf (c)}
are special cases of  Definition {\bf (a)}.  $C(p,q)=Q_{0}$ and $C^{a}(p,q)=Q_{n}$
are two arcs of $C$.

In continuous mathematics, the concept of the cross-over of two paths is called {\it transversal intersection}.
It means that one curve or path goes through (or penetrates) another curve.

\end{appendix}

%
%
%
%

\end{document}